\documentclass[a4paper,12pt]{article}
\usepackage[english]{babel}   %% загружает пакет многоязыковой вёрстки
\usepackage{indentfirst}
\usepackage[T2A]{fontenc}
\usepackage[utf8]{inputenc}

\renewcommand{\epsilon}{\ensuremath{\varepsilon}}
\renewcommand{\phi}{\ensuremath{\varphi}}
\renewcommand{\kappa}{\ensuremath{\varkappa}}
\renewcommand{\le}{\ensuremath{\leqslant}}

\renewcommand{\ge}{\ensuremath{\geqslant}}

\usepackage{amsmath,amsfonts,amssymb,amsthm,mathtools} % AMS
\usepackage{icomma} % "Умная" запятая: $0,2$ --- число, $0, 2$ --- перечисление

\usepackage{graphicx}  % Для вставки рисунков
\graphicspath{{images/}{images2/}}  % папки с картинками
\usepackage{geometry} % Задание полей страницы
\usepackage{setspace} % Интерлиньяж
\usepackage{lastpage} 
\usepackage{soul} % Модификаторы начертания
  \usepackage{hyperref}
\usepackage[usenames,dvipsnames,svgnames,table,rgb]{xcolor}
\begin{document}

\noindent  \begin{center} {\huge About connection of one matrix of composite\begin{spacing}{1.7}\end{spacing} numbers with Legendre’s conjecture} \end{center}

\textbf{}

\begin{center} {\large Garipov Ilshat Ilsurovich }\end{center}
\begin{spacing}{-1.7}\end{spacing}
\begin{center} \textit {\normalsize Russia, Republic of Tatarstan, Naberezhnye Chelny}\end{center}
\begin{spacing}{-1.7}\end{spacing}
\begin{center} \textit {\normalsize e-mail: mathsciencegaripovii@gmail.com.}\end{center}
\renewcommand{\abstractname}{Abstract}
\begin{abstract}
{\small  Scientific paper is devoted to establish connection of \textit{T}-matrix -- matrix of composite numbers $6h\pm 1$ in special view -- with Legendre's conjecture.}

\end{abstract}

\begin{center}
\textbf{\large List of symbols}
\end{center}
\begin{spacing}{0.8}\end{spacing}

${\rm \mathbb{N}} _{0} $ -- set of all natural numbers with zero.

${\rm \mathbb{N}} $ -- set of all natural numbers.

${\rm  \mathbb{Z}} $ -- set of all integers. 

${\rm \mathbb{P}} $ -- set of all prime numbers.

${\rm \mathbb{R}}$ -- set of all real numbers.

$T$ -- matrix comprising all defining and all not defining elements.

$\widetilde{T}$ -- set of all elements of T-matrix.

$D(b)$ -- $T$-matrix upper defining element of number $b$.

$d(b)$ -- $T$-matrix lower defining element of number $b$.

$D_{k} (b)$ -- upper defining element of number $b$ in $k$-row ($k>1$) of $T$-matrix.

$d_{k} (b)$ -- lower defining element of number $b$ in $k$-row ($k>1$) of $T$-matrix.

$W(b)$ -- $T$-matrix upper element of number $b$.

$w(b)$ -- $T$-matrix lower element of number $b$.

$W_{k} (b)$ -- upper element of number $b$ in $k$-row ($k>1$) of $T$-matrix.

$w_{k} (b)$ -- lower element of number $b$ in $k$-row ($k>1$) of $T$-matrix.

$D_{T} $ -- set of all defining elements of $T$-matrix.

$nD_{T} $ -- set of all not defining elements of $T$-matrix.

${\rm M} _{T} $ -- set of all leading elements of $T$-matrix.

$D_{T_{k} } $ -- set of all defining elements in $k$-row ($k>1$) of $T$-matrix.

$\pi (x)$ -- function counting the number of prime numbers less than or equal to $x\in {\rm \mathbb{R}}$.

$\pi _{{\rm M} _{T} } (x)$ -- function counting the number of $T$-matrix leading elements less than or equal to $x\in {\rm \mathbb{R}}$.

$\# _{k} (a)$ -- number of element $a$ in $k$-row of $T$-matrix.

${\rm H} _{(m-1)^{4} ,{\rm \; }m^{4} } $ -- «active» set for numbers $(m-1)^{4} ,{\rm \; }m^{4} $.\textbf{ }

$C_{(m-1)^{4} ,{\rm \; }m^{4} } $ -- «critical» element for numbers $(m-1)^{4} ,{\rm \; }m^{4} $.

$\nu_{k} (x)$ -- function counting the number of elements, less than or equal to $x\in {\rm \mathbb{R}}$, in $k$-row of $T$-matrix.

$\nu(x)$ -- function counting the number of naturals of the form $6h\pm 1$, less than or equal to $x\in {\rm \mathbb{R}}$.

$a\% b$ -- remainder after dividing $a\in {\rm \mathbb{N}}$ by $b\in {\rm \mathbb{N}} $.

$q_{m}$ -- number of prime numbers between $m^{2}$ and $(m+1)^{2} $.

\begin{center}
\textbf{\large Introduction 1. $T$-matrix}
\end{center}
\begin{spacing}{0.8}\end{spacing}

We construct a matrix $T\equiv \left(a(k;n)\right)_{\infty \times \infty } $, where $a(k;n)$ is a $T$-matrix element located in $k$-th row, $n$-th column and defined as follows:
\begin{spacing}{0.8}\end{spacing}
\[a(k;n)\equiv p(k)\cdot \left(5+2\cdot \left\lfloor \frac{n}{2} \right\rfloor +4\cdot \left\lfloor \frac{n-1}{2} \right\rfloor \right),\] 
\noindent where $p(k)$ is the $k$-th element of sequence $(p(k))_{k=1}^{\infty }$ of prime numbers: 
\begin{spacing}{0.5}\end{spacing}
\begin{equation} \label{(1)} 
p(k)\equiv p_{k+2} ,   
\end{equation}\begin{spacing}{1.2}\end{spacing}
\noindent where $p_{i}$ is the $i$-th prime number in sequence of all prime numbers (see [1]).

Let $(f(n))_{n=1}^{\infty } $ is a numerical sequence, where a common member $f(n)$ is defined as follows:
\begin{spacing}{0.8}\end{spacing}\[f(n)\equiv 3n+\frac{3-(-1)^{n} }{2} .\]\begin{spacing}{1.6}\end{spacing}

THEOREM 1.1. \begin{equation} \label{(2)} \left(\forall k,n\in {\rm \mathbb{N}} \right)\left(a(k;n)=p(k)\cdot f(n)\right).\end{equation}

DEFINITION 1.1. An element $a(k;n)$ of matrix $T$ is called defining if

1)  $a(k;n)$ is not divisible by 5;

2) $a(k;n)$ can be expressed as a product of some two prime numbers, that is
\begin{spacing}{0.5}\end{spacing}\begin{equation} \label{(3)} 
5\not |a(k;n){\rm \; \; }\wedge {\rm \; \; }(\exists p_{1} ,p_{2} \in {\rm \mathbb{P}} )(a(k;n)=p_{1} \cdot p_{2} ).  
\end{equation} \begin{spacing}{1.1}\end{spacing}

DEFINITION 1.2. An element $a(k;n)$ of matrix $T$ is called not defining if he does not satisfy condition \eqref{(3)}.

DEFINITION 1.3. An element $a(k;n)$ of matrix $T$ is called leading if
\begin{spacing}{0.5}\end{spacing}
\[a(k;n)=p^{2} (k).\]

DEFINITION 1.4. A $T$-matrix is called matrix comprising all defining and not defining elements.

LEMMA 1.2. $(f(n))_{n=1}^{\infty } $ is a sequence of all numbers of the form $6h\pm 1$:
\begin{spacing}{0.5}\end{spacing}
\[5;{\rm \; }7;{\rm \; }11;{\rm \; }13;{\rm \; }17;{\rm \; }19;{\rm \; }23;{\rm \; }25;{\rm \; }...{\rm \; };{\rm \; }6h-1;{\rm \; 6}h+{\rm 1;...}.\]
\begin{spacing}{1.1}\end{spacing}

PROPERTY 1.1. The sequence $(p^{2} (k))_{k=1}^{\infty }$ of $T$-matrix leading elements is ascending.

The simplest properties and basic theorems about elements of $T$-matrix are proved in [1].

\begin{center}
\textbf{\large 2. About a $T$-matrix upper defining element of real number}
\end{center}
\begin{spacing}{0.7}\end{spacing}

DEFINITION 2.1. A $T$-matrix defining element $D(b)$ is called an upper defining element of number $b\in {\rm \mathbb{R}}: b \ge 49$, if 
\begin{spacing}{0.5}\end{spacing}
\[D(b)=\min_{\scriptstyle \substack{a(k_{1}; n)\in D_{T} \\  a(k_{1}; n)>b \\  n\in {\rm \mathbb{N}}}}a(k_{1}; n), \]      
 \begin{spacing}{1.3}\end{spacing}

\noindent where $k_{1} $  is defined by condition
\begin{spacing}{0.3}\end{spacing}
\[ p^{2} (k_{1})=\max_{\scriptstyle \substack{p^{2} (k)\le b \\  k>1 }}{\rm \;}p^{2}(k).\]

DEFINITION 2.2. A $T$-matrix defining element $d(b)$ is called a lower defining element of number $b\in {\rm \mathbb{R}}: b>49$, if 
 \begin{spacing}{0.5}\end{spacing}\[d(b)=\max_{\scriptstyle \substack{ a(k_{2}; n)\in D_{T} \\  a(k_{2}; n)<b \\ n\in {\rm \mathbb{N}}}}a(k_{2}; n),\]

\noindent where $k_{2} $ is defined by condition
\begin{spacing}{0.3}\end{spacing}
\[ p^{2} (k_{2})= \max_{\scriptstyle \substack{ p^{2} (k)< b \\  k>1 }}{\rm \;}p^{2}(k).\]

DEFINITION 2.3. A $T$-matrix defining element $D_{k} (b)$ is called an upper defining element of number $b \in {\rm \mathbb{R}:}{\rm \; \;}p^{2} (k)\le b$, in $k$-row ($k>1$) of $T$-matrix if 
\begin{spacing}{0.5}\end{spacing}
\[D_{k} (b)=\min_{\scriptstyle \substack{ a(k; n)\in D_{T}  \\  a(k; n)>b \\  n\in {\rm \mathbb{N}}}}a(k; n). \]   
\begin{spacing}{1.3}\end{spacing}

DEFINITION 2.4. A $T$-matrix defining element $d_{k} (b)$ is called a lower defining element of number $b\in {\rm \mathbb{R}:}{\rm \; \;}p^{2} (k)<b$, in $k$-row ($k>1$) of $T$-matrix if
\begin{spacing}{0.5}\end{spacing}
\[d_{k} (b)=\max_{\scriptstyle \substack{ a(k; n) \in D_{T}  \\  a(k; n)<b \\ n \in {\rm \mathbb{N}}}} a(k; n).  \]  
\begin{spacing}{1.3}\end{spacing}

DEFINITION 2.5. A $T$-matrix element $W(b)$ is called an upper element of number $b\in {\rm \mathbb{R}}:$
\noindent $ b \ge 49$, if          
\begin{spacing}{0.1}\end{spacing}
\[W(b)=\min_{\scriptstyle \substack{ a(k_{1}; n) \in \widetilde{T}  \\  a(k_{1}; n)>b \\ n\in {\rm \mathbb{N}}}}a(k_{1}; n),  \]  
\begin{spacing}{1.3}\end{spacing}

\noindent where $k_{1} $ is defined by condition
\begin{spacing}{0.3}\end{spacing}
\[ p^{2} (k_{1})=\max_{\scriptstyle \substack{p^{2} (k)\le b \\ k>1}}{\rm \;}p^{2}(k).\]

DEFINITION 2.6. A $T$-matrix element $w(b)$ is called a lower element of number $b\in {\rm \mathbb{R}}:$

\noindent $b>49$, if                              
\begin{spacing}{0.5}\end{spacing}
\[w(b)=\max_{\scriptstyle \substack{ a(k_{2}; n)\in \widetilde{T} \\  a(k_{2}; n)<b \\ n\in {\rm \mathbb{N}}}}a(k_{2}; n),  \]  
\begin{spacing}{1.3}\end{spacing}

\noindent where $k_{2} $ is defined by condition
\begin{spacing}{0.3}\end{spacing}
\[ p^{2} (k_{2})= \max_{\scriptstyle \substack{ p^{2} (k)< b \\ k>1 }}{\rm \;}p^{2}(k).\]

DEFINITION 2.7. A $T$-matrix element $W_{k} (b)$ is called an upper element of number $b\in {\rm \mathbb{R}:}$

\noindent $p^{2} (k)\le b$, in $k$-row ($k>1$) of $T$-matrix if
\begin{spacing}{0.3}\end{spacing}
\[W_{k} (b)=\min_{\scriptstyle \substack{ a(k; n)\in \widetilde{T}  \\  a(k; n)>b \\ n\in {\rm \mathbb{N}}}}a(k; n).\]  

DEFINITION 2.8. A $T$-matrix element $w_{k} (b)$ is called a lower element of number $b\in {\rm \mathbb{R}:}$

\noindent $p^{2} (k)<b$, in $k$-row ($k>1$) of $T$-matrix if
\begin{spacing}{0.5}\end{spacing}
\[w_{k} (b)=\max_{\scriptstyle \substack{a(k; n) \in \widetilde{T} \\ a(k; n)<b \\ n \in {\rm \mathbb{N}}}}a(k; n).\]  
\begin{spacing}{1.3}\end{spacing}

LEMMA 2.1. 
\begin{spacing}{0.2}\end{spacing}
\begin{equation} \label{(4)} 
\left(\forall k>1\right)\left(\forall n>1\right)\left(a(k;n)\in D_{T} {\rm \; \; }\Leftrightarrow {\rm \; \; }f(n)\in {\rm \mathbb{P}} \backslash \{ 2;3;5\} \right).
\end{equation}
\begin{spacing}{1.2}\end{spacing}

PROOF. Choose any $k$-row ($k>1$) and any $n$-column ($n>1$) of $T$-matrix.

\textbf{Necessity. }Let $a(k;n)\in D_{T}$. By Definition 1.1, that means that
\begin{spacing}{0.5}\end{spacing}
\[5\not |a(k;n){\rm \; \; }\wedge {\rm \; \; }(\exists {\rm \; }p_{1} ,{\rm \; }p_{2} \in {\rm \mathbb{P}})(a(k;n)=p_{1} \cdot p_{2} ).\]
\begin{spacing}{1.2}\end{spacing}

Then by rule \eqref{(1)}, Lemma 1.2 and Theorem 1.1,
\begin{spacing}{0.5}\end{spacing}
\[p_{1} ={\rm \; }p(k)\in {\rm \mathbb{P}} \backslash \{ 2;3;5\} {\rm \; \; }\wedge {\rm \; \; }p_{2} =f(n)\in {\rm \mathbb{P}} \backslash \{ 2;3;5\}.\]
\begin{spacing}{1.5}\end{spacing}
It follows that $f(n)\in {\rm \mathbb{P}} \backslash \{ 2;3;5\} $. The necessity is proved.

\textbf{Sufficiency. }Let $f(n)\in {\rm \mathbb{P}} \backslash \{2;3;5\} $. It is clear that $p(k)\in {\rm \mathbb{P}} \backslash \{2;3;5\}$ for $k>1$.
\[{\rm \; \;}p(k),{\rm \;}f(n)\in {\rm \mathbb{P}} \backslash \{2;3;5\} {\rm \; }\Rightarrow{\rm \; }\]
\[\Rightarrow{\rm \; \; } p(k)\cdot f(n)\mathop{=}\limits^{\eqref{(2)}} a(k;n) {\rm \; \; }\wedge {\rm \; \;}p(k),{\rm \;}f(n)\in {\rm \mathbb{P}}{\rm \; \; }\wedge {\rm \; \; }5\not|a(k;n){\rm \; \; }\Rightarrow\]
\begin{spacing}{0.5}\end{spacing}
\[\Rightarrow{\rm \; \;}5\not |a(k;n){\rm \; \; }\wedge {\rm \; \; }(\exists {\rm \; }p_{1} ,p_{2} \in {\rm \mathbb{P}})(a(k;n)=p_{1} \cdot p_{2}){\rm \; \; } \mathop{\Leftrightarrow }\limits^{\eqref{(3)}}{\rm \; \;}a(k;n)\in D_{T}.\]

 The sufficiency is proved. Lemma 2.1 is proved. 

THEOREM 2.2 (about the «transition down» of $T$-matrix defining element). 
\begin{spacing}{0.5}\end{spacing}\[\left(\forall k;n\in {\rm \mathbb{N}} \right)\left(p^{2} (k)<a(k;n){\rm \; \; }\wedge {\rm \; \; }a(k;n)\in D_{T} {\rm \;}\Rightarrow \right. \] 
\begin{equation} \label{(5)} 
 {\rm \;}\left. \Rightarrow \left(\exists ! {\rm \; }j\in {\rm \mathbb{N}}\right)\left(k<j{\rm \;}\wedge {\rm \;}a(k;n)<p^{2} (j){\rm \;}\wedge {\rm \;}a(j;\# _{k} (p^{2} (k)))=a(k;n){\rm \;}\wedge {\rm \;}a(j;n)=p^{2} (j)\right)\right).
\end{equation}
\begin{spacing}{1.2}\end{spacing}

PROOF. \textbf{Existence.} It is established in [1].

\textbf{Uniqueness. }Suppose,
\begin{spacing}{0.5}\end{spacing}
\[(\forall k;n\in {\rm \mathbb{N}})(p^{2} (k)<a(k;n){\rm \; \; }\wedge {\rm \; \; }a(k;n)\in D_{T} {\rm \;}\Rightarrow  \] 
\[ \Rightarrow \left(\exists {\rm \; }j_{1}, j_{2} \in {\rm \mathbb{N}}\right)(j_{1} \ne j_{2}{\rm \;}\wedge {\rm \;}k<j_{1}{\rm \;}\wedge {\rm \;}k<j_{2}{\rm \;}\wedge{\rm \;}a(k;n)<p^{2} (j_{1}){\rm \;}\wedge {\rm \;} a(k;n)<p^{2} (j_{2}) {\rm \;}\wedge\]
\[\wedge{\rm \; \;}a(j_{1};\# _{k} (p^{2} (k)))=a(k;n){\rm \;}\wedge{\rm \;} a(j_{2};\# _{k} (p^{2} (k)))=a(k;n){\rm \;}\wedge \]
\[\wedge {\rm \; \;}a(j_{1};n)=p^{2} (j_{1}){\rm \;}\wedge{\rm \;}a(j_{2};n)=p^{2} (j_{2}))).\]
\[a(j_{1};\# _{k} (p^{2} (k)))=a(k;n){\rm \;}\wedge{\rm \;} a(j_{2};\# _{k} (p^{2} (k)))=a(k;n){\rm \;}\wedge{\rm \;} j_{1} \ne j_{2} {\rm \;}\Rightarrow {\rm \;}\] 
\[ \Rightarrow {\rm \;}a(j_{1};\# _{k} (p^{2} (k)))=a(j_{2};\# _{k} (p^{2} (k))){\rm \;}\wedge{\rm \;} j_{1} \ne j_{2} {\rm \; } \mathop{\Leftrightarrow }\limits^{\eqref{(2)}}\]
\[\Leftrightarrow{\rm \; }p(j_{1}) \cdot f(\# _{k} (p^{2} (k)))=p(j_{2}) \cdot f(\# _{k} (p^{2} (k))){\rm \;}\wedge{\rm \;} j_{1} \ne j_{2}{\rm \;}\Leftrightarrow\]
\[ {\rm \;}\Leftrightarrow{\rm \;} p(j_{1})=p(j_{2}){\rm \;}\wedge{\rm \;} j_{1} \ne j_{2}{\rm \;}\Leftrightarrow {\rm \;}j_{1}=j_{2}{\rm \;}\wedge{\rm \;} j_{1} \ne j_{2}.\]

As a result, a contradiction. The uniqueness is established. 

Theorem 2.2 is proved.

COROLLARY 2.3. Let $a(k;n)$, $a(j;n)$ are $T$-matrix elements from Theorem 2.2. Then,
\[\frac{a(k;n)}{p(k)} -p(k)=p(j)-\frac{a(k;n)}{p(j)} .\] 

PROOF. From algorithm №1 in [1], we get
\[\left(\exists h\in {\rm \mathbb{N}} \right)\left(p^{2} (k)+2h\cdot p(k)=a(k;n)\right){\rm \; \; }\Leftrightarrow {\rm \; \; }\left(\exists h\in {\rm \mathbb{N}} \right)\left(p(k)+2h=\frac{a(k;n)}{p(k)} \right).\] 
\[2h=\frac{a(k;n)}{p(k)} -p(k){\rm \; }\mathop{=}\limits^{ \eqref{(5)}} {\rm \; \; }\frac{a(j;{\rm \; }\# _{k} (p^{2} (k)))}{p(k)} -p(k){\rm \; \; }\mathop{=}\limits^{ \eqref{(2)}} {\rm \; \; }\frac{p(j)\cdot f(\# _{k} (p^{2} (k)))}{p(k)} -p(k)=\] 
\begin{spacing}{1.2}\end{spacing}
\[=\frac{p(j)\cdot p(k)}{p(k)} -f(\# _{k} (p^{2} (k))){\rm \; }\mathop{=}\limits^{ \eqref{(2)}} {\rm \; }p(j)-\frac{a(j;{\rm \; }\# _{k} (p^{2} (k)))}{p(j)} {\rm \; \; }\mathop{=}\limits^{ \eqref{(5)}} {\rm \; }p(j)-\frac{a(k;n)}{p(j)} .\] 
\begin{spacing}{1.7}\end{spacing}

Corollary 2.3 is proved.

CONCLUSION 2.1. Subject to the conditions of Theorem 2.2, there are the following equalities with some $h\in {\rm \mathbb{N}} $:
\begin{spacing}{0.4}\end{spacing}
\[1){\rm \;} p^{2} (k)+2h\cdot p(k)=a(k;n).\] 
\[2){\rm \;}a(k;n)+2h\cdot p(j)=p^{2} (j).\] 

Further, let $g_{k}^{-} \equiv p(k+1)-p(k)$.

LEMMA 2.4. 
\[(\forall k>1)(p^{2}(k)+g_{k}^{-} \cdot p(k)=D_{k}(p^{2}(k))).\]

PROOF. Choose any $k$-row ($k>1$) of $T$-matrix. Given Theorem 1.1, assume that 
\begin{spacing}{0.5}\end{spacing}
\[D_{k}(p^{2}(k))=p(k) \cdot f(n) \text{ with some } n>1.\]
\[D_{k}(p^{2}(k)) \in D_{T_{k}} {\rm \; }\mathop{\Rightarrow }\limits^{D_{T_{k}} \subset  D_{T}}{\rm \; }D_{k}(p^{2}(k)) \in D_{T}{\rm \; \; } \mathop{\Leftrightarrow }\limits^{\eqref{(4)}}{\rm \; \; }f(n)\in {\rm \mathbb{P}} \backslash \{ 2;3;5\}. \]
\begin{spacing}{1.2}\end{spacing}

It follows from Definition 2.3 that the defining elements of $T$-matrix don't exist between the elements $p^{2} (k)$ and $a(k;n)$. Therefore, $f(n)=p(k+1)$. Then,
\begin{spacing}{0.5}\end{spacing}\[D_{k}(p^{2}(k))= p(k)\cdot p(k+1)=p(k)\cdot (p(k)+g_{k}^{-} )=p^{2} (k)+g_{k}^{-} \cdot p(k).\] \begin{spacing}{1.2}\end{spacing}

Lemma 2.4 is proved.

COROLLARY 2.5. 
\begin{spacing}{0.3}\end{spacing}
\[(\forall k>1)(D_{k}(p^{2}(k))+g_{k}^{-} \cdot p(k+1)=p^{2} (k+1)).\]

PROOF. Choose any $k$-row ($k>1$) of $T$-matrix. Using Lemma 2.4, we get
\begin{spacing}{0.3}\end{spacing}
\[D_{k}(p^{2}(k))+g_{k}^{-} \cdot p(k+1)=p^{2}(k)+g_{k}^{-} \cdot p(k)+g_{k}^{-} \cdot p(k+1)=\]
\[=p^{2}(k)+g_{k}^{-} \cdot (p(k)+p(k+1))=p^{2}(k)+(p(k+1)-p(k))\cdot (p(k+1)+p(k))=\]
\[=p^{2}(k)+p^{2}(k+1)-p^{2}(k)=p^{2}(k+1).\]

Corollary 2.5 is proved. 

COMMENT. It follows from Definition 2.1 and Definition 2.3 that
\begin{spacing}{0.5}\end{spacing}
\[(\forall k>1)(D_{k}(p^{2}(k))=D(p^{2}(k))).\]
\begin{spacing}{1.2}\end{spacing}

THEOREM 2.6 (about the «transition up» of $T$-matrix defining element).
\begin{spacing}{0.5}\end{spacing}\[\left(\forall j;n\in {\rm \mathbb{N}} \right)\left(a(j;n)<p^{2} (j){\rm \; \; }\wedge {\rm \; \; }a(j;n)\in D_{T} {\rm \; \; }\Rightarrow \right. \] 
\[\left. \Rightarrow {\rm \;}\left(\exists ! {\rm \; }k\in {\rm \mathbb{N}}\right)\left(k<j{\rm \;}\wedge {\rm \;}p^{2} (k)<a(j;n){\rm \;}\wedge {\rm \;}a(k;{\rm \; }\# _{j} (p^{2} (j)))=a(j;n){\rm \;}\wedge {\rm \;}a(k;n)=p^{2} (k)\right)\right).\] \begin{spacing}{1.2}\end{spacing}

PROOF. \textbf{Existence. }Choose any defining element $a(j;n)$ that is smaller than the leading element $p^{2} (j)$ in $j$-row of $T$-matrix. Then it follows from Theorem 1.1 and Definition 1.1 that
\begin{spacing}{0.4}\end{spacing}
\[a(j;n)=p(j)\cdot f(n){\rm \; \;}\wedge {\rm \; \;}p(j),{\rm \; }f(n)\in {\rm \mathbb{P}} \backslash \{ 2;3;5\} .\] 
\begin{spacing}{1.2}\end{spacing} 

Therefore, the prime number $f(n)>5$ is an element of sequence $(p(k))_{k=1}^{\infty } $:
\begin{spacing}{0.5}\end{spacing}
\begin{equation} \label{(6)} 
\left(\exists {\rm \; }k\in {\rm \mathbb{N}} \right)\left(p(k)=f(n)\right).  
\end{equation} 
\begin{spacing}{1.2}\end{spacing}
$1) {\rm \;}a(j;n)<p^{2} (j){\rm \; \; }\mathop{\Leftrightarrow }\limits^{{\rm \eqref{(2)}}} {\rm \; \; }p(j)\cdot f(n)<p^{2} (j){\rm \; \; \; }\mathop{\Leftrightarrow }\limits^{{\rm \eqref{(6)}}} {\rm \; \; }p(j)\cdot p(k)<p^{2} (j){\rm \; \; }\Leftrightarrow$
\begin{spacing}{0.5}\end{spacing}
\[\Leftrightarrow{\rm \; \; }p(k)<p(j){\rm \; \; }\Leftrightarrow {\rm \; \; }k<j.\]

\textbf{}
\begin{spacing}{0.1}\end{spacing}
$2) {\rm \;}p(k)<p(j){\rm \; \;}\Leftrightarrow {\rm \; \;}p^{2} (k)<p(j)\cdot p(k){\rm \; \; \; }\mathop{\Leftrightarrow }\limits^{{\rm \eqref{(6)}}} {\rm \; \; }p^{2} (k)<p(j)\cdot f(n){\rm \; \; \; }\mathop{\Leftrightarrow }\limits^{{\rm \eqref{(2)}}} {\rm \; \; }p^{2} (k)<{\rm \; }a(j;n).$ 
\begin{spacing}{1.7}\end{spacing}
$3) {\rm \;}a(k;{\rm \; }\# _{j} (p^{2} (j)))\mathop{=}\limits^{\eqref{(2)}} p(k)\cdot f(\# _{j} (p^{2} (j)))=p(k)\cdot p(j)=p(j)\cdot p(k)\mathop{=}\limits^{\eqref{(6)}} p(j)\cdot f(n)\mathop{=}\limits^{\eqref{(2)}} a(j;n).$

\textbf{}
\begin{spacing}{0.5}\end{spacing}
$ 4){\rm \;}a(k;n)\mathop{=}\limits^{\eqref{(2)}} p(k)\cdot f(n){\rm \;}\mathop{=}\limits^{\eqref{(6)}} {\rm \; }p^{2} (k)$.

\textbf{Uniqueness. }Suppose,
\begin{spacing}{0.5}\end{spacing}
\[(\forall j;n\in {\rm \mathbb{N}})(a(j;n)<p^{2} (j){\rm \;}\wedge {\rm \;}a(j;n)\in D_{T} {\rm \; }\Rightarrow \] 
\[\Rightarrow (\exists {\rm \; }k_{1}, k_{2} \in {\rm \mathbb{N}})(k_{1} \ne k_{2}{\rm \;}\wedge {\rm \;} k_{1}<j{\rm \;}\wedge {\rm \;}k_{2}<j{\rm \;}\wedge {\rm \;}p^{2} (k_{1})<a(j;n){\rm \;}\wedge {\rm \;}p^{2} (k_{2})<a(j;n){\rm \;}\wedge\]
\[\wedge {\rm \; \;}a(k_{1};{\rm \; }\# _{j} (p^{2} (j)))=a(j;n){\rm \;}\wedge {\rm \;}a(k_{2};{\rm \; }\# _{j} (p^{2} (j)))=a(j;n){\rm \;}\wedge\]
\[ \wedge{\rm \; \;}a(k_{1};n)=p^{2} (k_{1}){\rm \;}\wedge {\rm \;}a(k_{2};n)=p^{2} (k_{2}))).\]
\[a(k_{1};{\rm \; }\# _{j} (p^{2} (j)))=a(j;n){\rm \;}\wedge {\rm \;}a(k_{2};{\rm \; }\# _{j} (p^{2} (j)))=a(j;n){\rm \;}\wedge {\rm \;}k_{1} \ne k_{2}{\rm \; }\Rightarrow\]
\[\Rightarrow{\rm \; }a(k_{1};{\rm \; }\# _{j} (p^{2} (j)))=a(k_{2};{\rm \; }\# _{j} (p^{2} (j))){\rm \;}\wedge{\rm \;} k_{1} \ne k_{2}{\rm \; } \mathop{\Leftrightarrow }\limits^{\eqref{(2)}}\]
\[\Leftrightarrow{\rm \; }p(k_{1}) \cdot f(\# _{j} (p^{2} (j)))=p(k_{2}) \cdot f(\# _{j} (p^{2} (j))){\rm \;}\wedge {\rm \;}k_{1} \ne k_{2}{\rm \; }\Leftrightarrow {\rm \;}\]
\[\Leftrightarrow{\rm \; }p(k_{1})=p(k_{2}) {\rm \;}\wedge {\rm \;}k_{1} \ne k_{2}{\rm \; }\Leftrightarrow{\rm \; }k_{1}=k_{2} {\rm \;}\wedge {\rm \;}k_{1} \ne k_{2}.\]

As a result, a contradiction. The uniqueness is established. 

Theorem 2.6 is proved.

LEGENDRE'S CONJECTURE.
\begin{spacing}{0.3}\end{spacing}\[\left(\forall m\in {\rm \mathbb{N}} \right)\left(\exists p\in {\rm \mathbb{P}} \right)\left(m^{2} <p<(m+1)^{2} \right).\]
\begin{spacing}{1.2}\end{spacing}

PROPOSITION 2.7. Any real number $x$ can be uniquely expressed as a sum of integer part (entire) and fractional part (mantissa) of number $x$:

\[x=\left\lfloor x\right\rfloor +\{ x\} .\] 
\begin{spacing}{0.8}\end{spacing}

PROPERTY 2.1 (property of number's entire). 
\begin{spacing}{0.3}\end{spacing}
\begin{equation} \label{(7)} (\forall x\in {\rm \mathbb{R}}){\rm \; }(x-1<\lfloor x \rfloor \le x).\end{equation}
\begin{spacing}{0.8}\end{spacing}

PROPERTY 2.2 (property of number's mantissa). 
\begin{spacing}{0.3}\end{spacing}
 \[(\forall x\in {\rm  \mathbb{R}}){\rm \; }({0}\le \{ x\} <1).\]
\begin{spacing}{0.8}\end{spacing}

PROPERTY 2.3. 
\begin{spacing}{0.3}\end{spacing}
\begin{equation} \label{(8)} \left(\forall n\in {\rm  \mathbb{N}} _{{\rm 0}} \right)\left(\forall x\in {\rm \mathbb{R}}\right)\left(n\cdot \left\lfloor x\right\rfloor \le \left\lfloor n\cdot x\right\rfloor \right) \text{ (see [2]).}\end{equation}

THEOREM 2.8 (fundamental theorem of arithmetic). Every positive integer except the number 1 can be represented in exactly one way apart from rearrangement as a product of one or more primes (see [3]).

PROPERTY 2.4. 
\begin{equation} \label{(9)}
\left(\forall x\in {\rm \mathbb{R}:\; \; }x\ge 0\right)\left(\nu(x)=\nu\left(\left\lfloor x\right\rfloor \right)\right).
\end{equation}

PROOF. Choose any real number $x\ge 0$. Then, using Proposition 2.7 and Property 2.2, we get 
\begin{spacing}{0.5}\end{spacing}
\[\nu(x)=\nu \left(\left\lfloor x\right\rfloor +\{ x\}\right)=\nu \left(\left\lfloor x\right\rfloor \right).\]

Property 2.4 is proved. 

THEOREM 2.9. 
\begin{spacing}{0.5}\end{spacing}
\begin{equation} \label{(10)}
\left(\forall k\in {\rm \mathbb{N}} \right)\left(\forall x\in {\rm \mathbb{R}:\; \;}x\ge 0\right)\left(\nu_{k} (x)=\nu\left(\frac{x}{p(k)} \right)\right).
\end{equation}
\begin{spacing}{1.5}\end{spacing}

PROOF. Fix any $k$-row of $T$-matrix and any real number $x\ge 0$. Select the elements 
\begin{spacing}{0.5}\end{spacing}
\[a(k;n)\in \widetilde{T}:{\rm \;}a(k;n)\le x,{\rm \; }n\in {\rm \mathbb{N}}.\]

Using Theorem 1.1 for each of them, we get the numbers $f(n)$:

\begin{spacing}{0.8}\end{spacing}
\begin{equation} \label{(11)} 
f(n)=\frac{a(k;n)}{p(k)} \le \frac{x}{p(k)} .     
\end{equation} 
\begin{spacing}{1.5}\end{spacing}

 By Lemma 1.2, the numbers $f(n)$ have the form $6h\pm 1$. Then, given \eqref{(11)}, the number $\nu(\frac{x}{p(k)})$ of such $f(n)$ is equal to the number $\nu_{k} (x)$ of elements, less than or equal to $x$, in $k$-row of $T$-matrix.

 Theorem 2.9 is proved. 

COROLLARY 2.10. 
 \begin{spacing}{0.2}\end{spacing}
 \[\left(\forall k\in {\rm \mathbb{N}} \right)\left(\forall x\in {\rm \mathbb{R}:\;\; }x\ge 0\right)\left(\nu_{k} (x)=\nu_{k} \left(\left\lfloor x\right\rfloor \right)\right).\]

PROOF. Fix any $k$-row of $T$-matrix and any real number $x\ge 0$. Then,
\begin{spacing}{0.9}\end{spacing}
\[{\rm \; \;}\nu_{k} (x){\rm \; }\mathop{=}\limits^{\eqref{(10)}} {\rm \;}\nu \left(\frac{x}{p(k)} \right){\rm \;}\mathop{=}\limits^{\eqref{(9)}} {\rm \;}\nu \left(\left\lfloor \frac{x}{p(k)} \right\rfloor \right){\rm \;}\mathop{=}\limits^{\eqref{(10)}} {\rm \;}\nu_{k} \left(p(k)\cdot \left\lfloor \frac{x}{p(k)} \right\rfloor \right){\rm \;}\mathop{\le }\limits^{{\rm \eqref{(8)}}}\]
\[ \le \nu_{k} \left(\left\lfloor p(k)\cdot \frac{x}{p(k)} \right\rfloor \right)=\nu_{k} \left(\left\lfloor x\right\rfloor \right).\] 
\begin{spacing}{1.5}\end{spacing}

It is clear that $\nu_{k}(x){\rm \;}\mathop{\ge }\limits^{{\rm \eqref{(7)}}} \nu_{k} \left(\left\lfloor x\right\rfloor \right)$. Thus,
 \begin{spacing}{0.5}\end{spacing}
\[\nu_{k} (x)\le \nu_{k} \left(\left\lfloor x\right\rfloor \right){\rm \;\;}\wedge{\rm \;\;} \nu_{k} (x)\ge \nu_{k} \left(\left\lfloor x\right\rfloor \right){\rm \;}\Rightarrow {\rm \;}\nu_{k} (x)=\nu_{k} \left(\left\lfloor x\right\rfloor \right).\]
 \begin{spacing}{1.2}\end{spacing}

 Corollary 2.10 is proved.

PROPOSITION 2.11. 
\begin{spacing}{0.8}\end{spacing}\begin{equation} \label{(12)} 
\left(\forall m\in {\rm \mathbb{N}} \right)\left(\nu(m)=\left\lfloor \frac{m+2}{3} \right\rfloor -\left\lfloor \frac{m\% 6}{4} \right\rfloor +\left\lfloor \frac{m\% 6}{5} \right\rfloor -1\right).
\end{equation}  \begin{spacing}{1.5}\end{spacing}

PROOF. The exact formula of space complexity $C(m)$ of $T$-matrix - based algorithm (algorithm №1) for finding all the prime numbers less than or equal to a given natural number $m\ge 5$ was obtained in [1]:
\begin{spacing}{0.8}\end{spacing}
\[\left(\forall m\in {\rm \mathbb{N}}:{\rm \;}m\ge 5\right)\left(C(m)=\left\lfloor \frac{m+2}{3} \right\rfloor -\left\lfloor \frac{m\% 6}{4} \right\rfloor +\left\lfloor \frac{m\% 6}{5} \right\rfloor \right).\]
\begin{spacing}{1.5}\end{spacing}

This formula also takes into account the number 0 for correct numbering of natural numbers of the form $6h\pm 1$ less than or equal to $m$. Discarding the number 0, we get
\begin{spacing}{0.5}\end{spacing}
\begin{equation} \label{(13)} 
(\forall m\in {\rm \mathbb{N}}:{\rm \;}m\ge 5)(\nu(m)=C(m)-1). 
\end{equation} 
\begin{spacing}{-0.5}\end{spacing}
\[\left(\forall m\in \{ 1;2;3;4\} \right)\left(\left\lfloor \frac{m+2}{3} \right\rfloor -\left\lfloor \frac{m\% 6}{4} \right\rfloor +\left\lfloor \frac{m\% 6}{5} \right\rfloor -1=0{\rm \; \; }\wedge {\rm \; \; }\nu(m)=0\right).\] 

 Proposition 2.11 is proved.

Further, we will present a method №1 which allows to find a $T$-matrix upper defining element $D(m^{4})$ of number $m^{4}$ ($m\in {\rm \mathbb{N}} :{\rm \;}m\ge 3$) by invoking that between $(m-1)^{2} $ and $m^{2} $ there is the prime number. Also note that a prime number between $m^{2} $ and $(m+1)^{2} $ is the intermediate result of method №1. The important comments within $T$-matrix are also given in method №1.  

\textbf{Description of method №1. Input:} $m\in {\rm \mathbb{N}}:{\rm \;}m\ge 3$. 

\textbf{Step 1. }Using Proposition 2.11, compute a number $\overline{n}$ of numbers of the form $6h\pm 1$ less than or equal to $m^{2}$:
\begin{spacing}{0.8}\end{spacing}
\[\overline{n}\equiv \nu(m^{2})=\left\lfloor \frac{m^{2}+2}{3} \right\rfloor -\left\lfloor \frac{m^{2}\% 6}{4} \right\rfloor +\left\lfloor \frac{m^{2}\% 6}{5} \right\rfloor-1.\] 
\begin{spacing}{1.5}\end{spacing}

\textbf{Step 2. }Find a prime number $p(k_{1})$:
\begin{spacing}{0.5}\end{spacing}
\begin{equation} \label{(14)} 
p(k_{1})=\max_{\scriptstyle \substack {(m-1)^{2}<p(k)<m^{2} \\ k>1}} {\rm \;}p(k).
\end{equation} 
\begin{spacing}{1.2}\end{spacing}

It follows that $(m-1)^{2} <p(k_{1} )<m^{2}$. These inequalities are strict, since $(m-1)^{2} ,{\rm \; }m^{2} \notin {\rm \mathbb{P}} $.

Within $T$-matrix,
\begin{spacing}{0.5}\end{spacing}
\[(m-1)^{2}<p(k_{1} )<m^{2} {\rm \; \; }\Leftrightarrow {\rm \; \; }(m-1)^{4}<p^{2} (k_{1} )<m^{4}, \text{where } p^{2} (k_{1} )\in {\rm M} _{T}. \]

Test the primality of numbers $f(\overline{n}-i)$ that lie between $(m-1)^{2} $ and $m^{2}$ starting at $i=0$ with a step 1 until a number $\Delta \overline{n}\in {\rm \mathbb{N}}_{0}$ is found: 
\begin{spacing}{0.5}\end{spacing}
\[f(\overline{n}-\Delta \overline{n})\in {\rm \mathbb{P}} \backslash \{ 2;3;5\}.\]

For numbers $f(\overline{n}-i){\rm \;}(i=\overline{0;\Delta \overline{n}})$ we use the Lenstra-Pomerance primality test (modification of polynomial-time Agrawal–Kayal–Saxena (AKS) primality test, see [4], [5]). Note that by Lemma 1.2, the numbers $f(\overline{n}-i)$ have the form $6h\pm 1$. Thus, 
\begin{spacing}{0.3}\end{spacing}\[p(k_{1})=f(\overline{n}-\Delta \overline{n}).\] \begin{spacing}{1}\end{spacing}

\textbf{Step 3. }Compute a number $n_{0}$ of $T$-matrix upper element $W(m^{4})$ of $m^{4}$ (in $k_{1}$-row of $T$-matrix).

Let's say that $a(k_{1} ;n_{0})\equiv W(m^{4} )$. Then by Definition 2.5,
\begin{spacing}{0.5}\end{spacing}\[a(k_{1} ;n_{0} )=\min_{\scriptstyle \substack { a(k_{1}; n)\in \widetilde{T}  \\  a(k_{1}; n)>m^{4} \\ n\in {\rm \mathbb{N}}}}a(k_{1}; n).\]

It follows that 
\begin{spacing}{0.1}\end{spacing}
\[\nu_{k_{1} } (m^{4} )+1=\nu_{k_{1} } (a(k_{1} ;n_{0} )){\rm \; \; }\wedge {\rm \; \; }n_{0} =\nu_{k_{1} } (a(k_{1} ;n_{0} )).\] 

Then,
\[n_{0} =\nu_{k_{1} } (m^{4} )+1.\]
\begin{spacing}{0.8}\end{spacing}
\[\nu_{k_{1} } (m^{4} ){\rm \;}\mathop{=}\limits^{\eqref{(10)}} {\rm \;}\nu\left(\frac{m^{4} }{p(k_{1} )} \right){\rm \;}\mathop{=}\limits^{\eqref{(9)}} {\rm \; \; }\nu\left(\left\lfloor \frac{m^{4} }{p(k_{1} )} \right\rfloor \right) {\rm \; \; }\mathop{\Rightarrow }{\rm \; \; }\nu\left(\left\lfloor \frac{m^{4} }{p(k_{1} )} \right\rfloor \right)=n_{0}-1.\] 

\[\left\lfloor \frac{m^{4} }{p(k_{1} )} \right\rfloor {\rm \;}\mathop{>}\limits^{\eqref{(7)}} {\rm \;}\frac{m^{4} }{p(k_{1} )} -1{\rm \;}\mathop{>}\limits^{\eqref{(14)}} {\rm \;}\frac{m^{4} }{m^{2} } -1=m^{2} -1\ge 3^{2} -1=8>5{\rm \; \;}\mathop{\Rightarrow }\limits^{\eqref{(13)}}\] 

\[\mathop{\Rightarrow }{\rm \; \;}\nu\left(\left\lfloor \frac{m^{4}}{p(k_{1})} \right\rfloor \right)=C\left(\left\lfloor \frac{m^{4}}{p(k_{1} )} \right\rfloor \right)-1{\rm \; \;}\mathop{\Leftrightarrow }{\rm \; \;}n_{0}-1=C\left(\left\lfloor \frac{m^{4} }{p(k_{1} )} \right\rfloor \right)-1{\rm \; \;}\mathop{\Leftrightarrow }\] 
\begin{spacing}{1.5}\end{spacing}
\[\mathop{\Leftrightarrow }{\rm \; \;}n_{0}=C\left(\left\lfloor \frac{m^{4} }{p(k_{1} )} \right\rfloor \right).\]

\textbf{Step 4. }Compute a $T$-matrix element $D(m^{4})$, and $D(m^{4})<(m+1)^{4}$.

\textbf{Test way to calculate the element $D(m^{4})$. }Given Definition 1.1, finding the element $D(m^{4})$ is reduced to test the primality of numbers $f(n_{0}+i)=\frac{a(k_{1}, n_{0}+i)}{p(k_{1})}$ starting at $i=0$ with a step 1 until a number $\Delta n_{0} \in {\rm \mathbb{N}}_{0}$ is found:
\begin{spacing}{0.5}\end{spacing}
\[f(n_{0} +\Delta n_{0} )\in {\rm \mathbb{P}} \backslash \{ 2;3;5\}.\]

In turn, by Lemma 2.1,
\begin{spacing}{0.5}\end{spacing}
\[f(n_{0} +\Delta n_{0} )\in {\rm \mathbb{P}} \backslash \{ 2;3;5\} {\rm \; \; \; }\Leftrightarrow {\rm \; \; \; }a(k_{1} ,n_{0} +\Delta n_{0} )\in D_{T} .\] 

For numbers $f(n_{0}+i) {\rm \;} (i=\overline{0;\Delta n_{0}})$ use the Lenstra-Pomerance primality test. Using Theorem 1.1, compute the element $D(m^{4})$:
\begin{spacing}{0.5}\end{spacing}
\[D(m^{4} )\equiv a(k_{1} ,n_{0} +\Delta n_{0} )=p(k_{1} )\cdot f(n_{0} +\Delta n_{0} ).\] 
\begin{spacing}{1.2}\end{spacing}

We expect that the received prime number $p(j)\equiv f(n_{0} +\Delta n_{0})$ lies between $m^{2}$ and $(m+1)^{2}$.
 
\textbf{Output: }$T$-matrix upper defining element $D(m^{4})$ of number $m^{4}$.

Next, we will use the Lagarias-Odlyzko analytical method (about the method see [5], [6]) to know which $k_{1}$, $j$-rows of $T$-matrix contain the element $D(m^{4})$.

The calculation of numbers $k_{1}$, $j$-rows of $T$-matrix is reduced to calculation of numbers $\overline{k}_{1} $ of $p(k_{1})$, 
$\overline {j}$ of $p(j)$ in sequence of all prime numbers respectively.

Clear that
\begin{spacing}{0.3}\end{spacing}
\[\overline{k}_{1}=\pi (p(k_{1})) {\rm \; \; } \wedge {\rm \; \; } \overline { j}=\pi(p(j)).\]

 Using rule (1), 
 \begin{spacing}{0.3}\end{spacing}
 \[k_{1}=\overline{k}_{1}-2  {\rm \; \; } \wedge {\rm \; \; } j=\overline { j}-2.\]

Consider now the inverse problem. Suppose that Step 1 and Step 2 of method №1 are done, but we already know which elements are defining, and which elements are not defining in $T$-matrix. In this case, a prime number $p(j)$ that lies between $m^{2}$ and $(m+1)^{2}$ can be found in two ways.
\begin{spacing}{1.7}\end{spacing}
\textbf{Way №1 (hard).}
\begin{spacing}{0.8}\end{spacing}
\begin{equation} \label{(15)}
p(j)=h+\sqrt{h^{2}+D(m^{4})},\text{ where } h=\frac{D(m^{4})-p^{2} (k_{1})}{2\cdot p(k_{1})}.                             
\end{equation}

PROOF. Compute a number $h\in {\rm \mathbb{N}}$ on the basis of Conclusion 2.1, 1):
\[p^{2} (k_{1})+2h\cdot p(k_{1})=D(m^{4} ){\rm \; \; }\Leftrightarrow {\rm \; \; }2h\cdot p(k_{1})=D(m^{4} )-p^{2} (k_{1}){\rm \; \; }\Leftrightarrow {\rm \; \; }h=\frac{D(m^{4} )-p^{2} (k_{1} )}{2\cdot p(k_{1} )} .\] 

On the basis of Conclusion 2.1, 2),
\begin{spacing}{0.5}\end{spacing}
\begin{equation} \label{(16)} 
D(m^{4} )+2h\cdot p(j)=p^{2} (j){\rm \; \; }\Leftrightarrow {\rm \; \; }p^{2} (j)-2h\cdot p(j)-D(m^{4} )=0.  
\end{equation} 
\begin{spacing}{1.2}\end{spacing}

Solve the reduced quadratic equation \eqref{(16)} in the unknown $p(j)$.
\[p_{1,2} (j)=\frac{2h\pm \sqrt{(-2h)^{2} -4\cdot (-D(m^{4} ))} }{2} =\frac{2h\pm \sqrt{4h^{2} +4\cdot D(m^{4} )} }{2} =\] 
\[=\frac{2h\pm 2\cdot \sqrt{h^{2}+D(m^{4} )}}{2}=h\pm \sqrt{h^{2} +D(m^{4} )} ,{\rm \; \; }h^{2} +D(m^{4} )>0.\] 

In turn, $h-\sqrt{h^{2}+D(m^{4})}<0$. Really,
\[D(m^{4} )>0{\rm \; \; }\Leftrightarrow {\rm \; \; }D(m^{4} )+h^{2} >h^{2} {\rm \; \; }\Leftrightarrow {\rm \; \; }\sqrt{D(m^{4} )+h^{2} } >\sqrt{h^{2} } {\rm \; \; }\Leftrightarrow\]
\[\Leftrightarrow{\rm \; \; }\sqrt{D(m^{4} )+h^{2} }>\left|h\right|.\] 
\begin{spacing}{0.5}\end{spacing}
\[h>0 {\rm \; \; }\Rightarrow {\rm \; \; }\sqrt{D(m^{4} )+h^{2} } >h{\rm \; \; }\Leftrightarrow {\rm \; \; }h-\sqrt{h^{2} +D(m^{4} )} <0.\] 

Thus, given $p(j)\ge 5$, the root $p(j)=h-\sqrt{h^{2} +D(m^{4} )} $ of equation \eqref{(16)} is not considered. 

Show that $h+\sqrt{h^{2} +D(m^{4} )}>5$. It follows from Definition 2.1 that $m^{4}<D(m^{4})$.
\begin{spacing}{-0.2}\end{spacing}\[h+\sqrt{h^{2} +D(m^{4} )} >h+\sqrt{h^{2} +m^{4} } \ge h+\sqrt{h^{2} +3^{4} } >h+\sqrt{3^{4} } =h+9>5.\] \begin{spacing}{0.9}\end{spacing}

In result, $p(j)=h+\sqrt{h^{2}+D(m^{4})}$ is an appropriate root of equation \eqref{(16)}. 

The correctness of formula \eqref{(15)} is proved.

\textbf{Way №2 (easy).}
\begin{spacing}{0.5}\end{spacing}
\begin{equation} \label{(17)}
p(j)=\frac{D(m^{4} )}{p(k_{1} )}.
\end{equation} 

Formula \eqref{(17)} can be obtained from the proof of Theorem 2.2 (see [1]).

Verify an equation of the prime numbers $p(j)$ found in 2 ways.
\begin{spacing}{1.5}\end{spacing}
\textbf{Verification. }$p^{2} (k_{1} )+2h\cdot p(k_{1} )=D(m^{4} ){\rm \; }\Leftrightarrow {\rm \; }(p^{2} (k_{1} )+2h\cdot p(k_{1} ))\cdot D(m^{4} )=D^{2} (m^{4} ) \Leftrightarrow $
\begin{spacing}{0.2}\end{spacing}
\[\Leftrightarrow {\rm \;}p^{2} (k_{1} )\cdot D(m^{4} )+2h\cdot p(k_{1} )\cdot D(m^{4} )=D^{2} (m^{4} ){\rm \; \; }\Leftrightarrow {\rm \; \; }D(m^{4} )+2h\cdot \frac{D(m^{4} )}{p(k_{1} )} =\frac{D^{2} (m^{4} )}{p^{2} (k_{1} )} {\rm \;}\Leftrightarrow\]
\[\Leftrightarrow {\rm \;}D(m^{4} )=\frac{D^{2} (m^{4} )}{p^{2} (k_{1} )} -2h\cdot \frac{D(m^{4} )}{p(k_{1} )} {\rm \; \; }\Leftrightarrow {\rm \; \; }D(m^{4} )+h^{2} =\frac{D^{2} (m^{4} )}{p^{2} (k_{1} )} -2h\cdot \frac{D(m^{4} )}{p(k_{1} )} +h^{2} {\rm \; \; }\Leftrightarrow \] 
\begin{spacing}{0.3}\end{spacing}
\[\Leftrightarrow {\rm \; \; }D(m^{4})+h^{2} =\left(\frac{D(m^{4})}{p(k_{1} )}-h\right)^{2} {\rm \; }\Leftrightarrow {\rm \; \; }
\sqrt{D(m^{4})+h^{2}}=\sqrt{\left(\frac{D(m^{4})}{p(k_{1} )}-h\right)^{2}} {\rm \; \; }\Leftrightarrow {\rm \; \; }\]
\begin{spacing}{1.2}\end{spacing}
\[\Leftrightarrow {\rm \; \; }\sqrt{D(m^{4} )+h^{2} } =\left|\frac{D(m^{4} )}{p(k_{1} )}-h\right|.\] 
\begin{spacing}{0.6}\end{spacing}
\[{\rm \; \; }p^{2} (k_{1} )+2h\cdot p(k_{1} )=D(m^{4} ){\rm \;}\Leftrightarrow {\rm \;}p(k_{1} )+2h=\frac{D(m^{4} )}{p(k_{1} )} {\rm \;}\Leftrightarrow{\rm \; }p(k_{1} )+h=\frac{D(m^{4} )}{p(k_{1} )}-h.\]
\[p(k_{1})+h>0{\rm \; \;}\Rightarrow{\rm \; \; } \frac{D(m^{4})}{p(k_{1})}-h>0{\rm \; \; }\Rightarrow\]
\begin{spacing}{1.2}\end{spacing}
\[\Rightarrow {\rm \; \; } \sqrt{D(m^{4} )+h^{2}}=\frac{D(m^{4} )}{p(k_{1} )}-h{\rm \; \; }\Leftrightarrow {\rm \; \; }h+\sqrt{h^{2} +D(m^{4})} =\frac{D(m^{4})}{p(k_{1})}.\] 

We will put forward a following conjecture.

CONJECTURE 2.1. 
\[\left(\forall m\in {\rm \mathbb{N}}:{\rm \;} m\ge 3\right)\left(m^{2} <\frac{D(m^{4} )}{p(k_{1} )} <(m+1)^{2} \right),\]

\noindent where $k_{1}$ is defined by condition \eqref{(14)}.

If Conjecture 2.1 is true, then for both ways:
\begin{spacing}{0.5}\end{spacing}
\[m^{4} <D(m^{4} )<p^{2} (j)<(m+1)^{4}.\]

Therefore, $m^{2}<p(j)<(m+1)^{2}$.

THEOREM 2.12. 
\begin{spacing}{0.2}\end{spacing}
\[\left(\forall m\in {\rm \mathbb{N}}:{\rm \;} m\ge 3\right)(p(k_{1} )\not {\rm |}{\rm \; }m^{4}),\]

\noindent where $k_{1}$ is defined by condition \eqref{(14)}.

PROOF. Suppose otherwise: $\left(\exists m\in {\rm \mathbb{N}}:{\rm \;} m\ge 3\right)(p(k_{1} ){\rm \; |}m^{4})$. 

Given Theorem 2.8, we present a number $m$ as
\begin{spacing}{0.8}\end{spacing}
\[m=\prod _{i=1}^{w}p_{\alpha _{i} }^{\beta _{i} }, \text{ where } w,{\rm \;}\alpha _{i} ,{\rm \;}\beta _{i} \in {\rm \mathbb{N}} ;{\rm \; }i=\overline{1;w}.\]
\begin{spacing}{1.2}\end{spacing}
\[m^{4} =\left(\prod _{i=1}^{w}p_{\alpha _{i} }^{\beta _{i} }\right)^{4} =\left(\prod _{i=1}^{w}\prod _{j=1}^{\beta _{i}}p_{\alpha _{i}} \right)^{4} =\prod _{i=1}^{w}\prod _{j=1}^{\beta _{i} }p_{\alpha _{i}}^{4} {\rm \; \; }\mathop{\Rightarrow }\limits^{p(k_{1} ){\rm \; |\; }m^{4} }\]
\begin{spacing}{1.0}\end{spacing}
\[\Rightarrow {\rm \; \; (}\exists i\in {\rm \mathbb{N}} :1\le i\le w)(p_{\alpha _{i} } =p(k_{1} )){\rm \; \; }\Leftrightarrow {\rm \; \; }p(k_{1} ){\rm \; |\; }m.\] 

Without loss of generality, assume that $m=p(k_{1} )\cdot m_{1} ,{\rm \; }m_{1} \in {\rm \mathbb{N}}$.

 It is known that $(\forall i\in {\rm \mathbb{N}} )(p_{i+1}<2\cdot p_{i})$ (see [7]). Therefore, using \eqref{(1)}, we get
 \begin{spacing}{0.5}\end{spacing}
\[p(k_{1} +1)<2 \cdot p(k_{1})<p^{2} (k_{1}) \cdot m_{1}^{2} =m^{2} .\] 
\begin{spacing}{1.2}\end{spacing}

 As a result, a contradiction to the maximality of the prime number $p(k_{1})$ from \eqref{(14)}. 

 Theorem 2.12 is proved.
 
PROPOSITION 2.13. 
\begin{spacing}{0.2}\end{spacing}
\begin{equation} \label{(18)}
\nu_{k_{1} } (D(m^{4} ))-\nu_{k_{1} } (m^{4} )=O(m),
\end{equation}

\noindent where $k_{1}$ is defined by condition \eqref{(14)}.

PROOF.Given Conjecture 2.1, let's say that 
 \begin{spacing}{0.5}\end{spacing}
 \[D(m^{4} )<p^{2} (j)<(m+1)^{4}.\]
 \begin{spacing}{-0.3}\end{spacing}
\[\nu_{k_{1} } (D(m^{4} ))-\nu_{k_{1} } (m^{4} ){\rm \; \; }\mathop{=}\limits^{\eqref{(10)}} {\rm \; \; }\nu\left(\frac{D(m^{4} )}{p(k_{1} )} \right)-\nu\left(\frac{m^{4} }{p(k_{1} )} \right){\rm \;}\mathop{=}\limits^{\eqref{(17)}, \eqref{(9)}}{\rm \;}\nu(p(j))-\nu\left(\left\lfloor \frac{m^{4} }{p(k_{1} )} \right\rfloor \right).\] 
\[\left(\forall n\in {\rm \mathbb{N}} \right)\left(-\left\lfloor \frac{n\% 6}{4} \right\rfloor +\left\lfloor \frac{n\% 6}{5} \right\rfloor \in \{ -1;0\} \right).\]
\begin{spacing}{2}\end{spacing}
Then,
\begin{spacing}{1.0}\end{spacing}
\[\nu(p(j))-\nu\left(\left\lfloor \frac{m^{4} }{p(k_{1} )} \right\rfloor \right){\rm \;}\mathop{=}\limits^{\eqref{(12)}} {\rm \;}\left\lfloor \frac{p(j)+2}{3} \right\rfloor -\left\lfloor \frac{p(j)\% 6}{4} \right\rfloor +\left\lfloor \frac{p(j)\% 6}{5} \right\rfloor -1-\] 
\begin{spacing}{1.2}\end{spacing}
\[-\left({\rm \; }\left\lfloor \frac{\left\lfloor m^{4} /p(k_{1} )\right\rfloor +2}{3} \right\rfloor -\left\lfloor \frac{\left\lfloor m^{4} /p(k_{1} )\right\rfloor \% 6}{4} \right\rfloor +\left\lfloor \frac{\left\lfloor m^{4} /p(k_{1})\right\rfloor \% 6}{5} \right\rfloor -1\right) \le\] 
\begin{spacing}{0.2}\end{spacing}
\[\le \left\lfloor \frac{p(j)+2}{3} \right\rfloor -1-\left\lfloor \frac{\left\lfloor m^{4} /p(k_{1} )\right\rfloor +2}{3} \right\rfloor -\left(\left\lfloor \frac{\left\lfloor m^{4} /p(k_{1} )\right\rfloor \% 6}{5} \right\rfloor -\left\lfloor \frac{\left\lfloor m^{4} /p(k_{1} )\right\rfloor \% 6}{4} \right\rfloor \right)+1=\] 
\[=\left\lfloor \frac{p(j)+2}{3} \right\rfloor -\left\lfloor \frac{\left\lfloor m^{4} /p(k_{1} )\right\rfloor +2}{3} \right\rfloor -\left(\left\lfloor \frac{\left\lfloor m^{4} /p(k_{1} )\right\rfloor\% 6}{5} \right\rfloor -\left\lfloor \frac{\left\lfloor m^{4} /p(k_{1} )\right\rfloor \% 6}{4} \right\rfloor \right) \le \] 
\begin{spacing}{0.6}\end{spacing}
\[\le \left\lfloor \frac{p(j)+2}{3} \right\rfloor -\left\lfloor \frac{\left\lfloor m^{4} /p(k_{1} )\right\rfloor +2}{3} \right\rfloor +1{\rm \;}\mathop{\le }\limits^{\eqref{(7)}} {\rm \;}\frac{p(j)+2}{3} -\left\lfloor \frac{\left\lfloor m^{4} /p(k_{1} )\right\rfloor +2}{3} \right\rfloor +1.\] 
\begin{spacing}{1.0}\end{spacing}
\[\left\lfloor \frac{\left\lfloor m^{4} /p(k_{1} )\right\rfloor +2}{3} \right\rfloor {\rm \;}\mathop{>}\limits^{\eqref{(7)}} {\rm \;}\left(\frac{1}{3} \cdot \left\lfloor \frac{m^{4} }{p(k_{1} )} \right\rfloor +\frac{2}{3} \right)-1=\frac{1}{3} \cdot \left\lfloor \frac{m^{4} }{p(k_{1} )} \right\rfloor -\frac{1}{3} {\rm \; \; }\Rightarrow \] 
\begin{spacing}{1.4}\end{spacing}
\[\Rightarrow {\rm \;}-\left\lfloor \frac{\left\lfloor m^{4} /p(k_{1} )\right\rfloor +2}{3} \right\rfloor <\frac{1}{3} -\frac{1}{3} \cdot \left\lfloor \frac{m^{4} }{p(k_{1} )} \right\rfloor {\rm \; \; }\Rightarrow{\rm \; \; }\]
\begin{spacing}{1.2}\end{spacing}
\[\Rightarrow {\rm \; \; }\nu(p(j))-\nu\left(\left\lfloor \frac{m^{4} }{p(k_{1} )} \right\rfloor \right)<\frac{p(j)}{3} +\frac{2}{3} +\left(\frac{1}{3} -\frac{1}{3} \cdot \left\lfloor \frac{m^{4} }{p(k_{1} )} \right\rfloor \right)+1{\rm \; \; }\Leftrightarrow {\rm \; \; }\]
\begin{spacing}{1.2}\end{spacing}
\begin{equation} \label{(19)} 
\Leftrightarrow {\rm \; \; }\nu(p(j))-\nu\left(\left\lfloor \frac{m^{4} }{p(k_{1} )} \right\rfloor \right)<\frac{p(j)}{3} -\frac{1}{3} \cdot \left\lfloor \frac{m^{4} }{p(k_{1} )} \right\rfloor +2.               
\end{equation} 
\[\left\lfloor \frac{m^{4} }{p(k_{1})} \right\rfloor {\rm \;}\mathop{>}\limits^{{\rm \eqref{(7)}}} {\rm \;}\frac{m^{4} }{p(k_{1} )} -1{\rm \; \;}\Leftrightarrow {\rm \; \;}-\left\lfloor \frac{m^{4} }{p(k_{1} )} \right\rfloor<1-\frac{m^{4} }{p(k_{1} )} {\rm \; \;}\Rightarrow {\rm \; \;}\nu(p(j))-\nu\left(\left\lfloor \frac{m^{4} }{p(k_{1} )} \right\rfloor \right)<\] 
\begin{spacing}{0.6}\end{spacing}
\[<\frac{p(j)}{3} +\frac{1}{3} \cdot \left(1-\frac{m^{4} }{p(k_{1})} \right)+2=\frac{p(j)}{3} +\frac{1}{3} -\frac{1}{3} \cdot \frac{m^{4} }{p(k_{1} )} +2<\frac{1}{3} \cdot (m+1)^{2} -\frac{1}{3} \cdot \frac{m^{4} }{p(k_{1} )} +\frac{7}{3} {\rm \;}.\] 
\[\frac{m^{4} }{p(k_{1} )} {\rm \; \; }\mathop{>}\limits^{\eqref{(14)}} {\rm \; \; }\frac{m^{4} }{m^{2} } =m^{2} {\rm \; \; }\Leftrightarrow {\rm \; \; }-\frac{m^{4} }{p(k_{1} )} <-{\rm \; }m^{2} {\rm \; \; }\Rightarrow {\rm \; \; }\nu(p(j))-\nu\left(\left\lfloor \frac{m^{4} }{p(k_{1} )} \right\rfloor \right)<\] 
\begin{spacing}{1.2}\end{spacing}
\[<\frac{1}{3} \cdot (m+1)^{2} -\frac{1}{3} \cdot m^{2} +\frac{7}{3} =\frac{1}{3} \cdot ((m+1)^{2} -m^{2} )+\frac{7}{3}=\]
\begin{spacing}{1.0}\end{spacing}
\[=\frac{1}{3} \cdot (m+1+m)\cdot (m+1-m)+\frac{7}{3}=\frac{1}{3} \cdot (2m+1)+\frac{7}{3}=\frac{2m}{3} +\frac{8}{3}<\]
\begin{spacing}{0.1}\end{spacing}
\begin{equation} \label{(20)} 
<\frac{2m}{3} +\frac{3m}{3} =\frac{5}{3} \cdot m{\rm \; \; }(m\ge 3) {\rm \; \; }\Rightarrow{\rm \; \; }\nu_{k_{1} } (D(m^{4} ))-\nu_{k_{1} } (m^{4} )<\frac{5}{3} \cdot m{\rm \; \; }(m\ge 3).                        
\end{equation} 

It follows from (20) that $\nu_{k_{1}}(D(m^{4}))-\nu_{k_{1}}(m^{4})=O(m)$. 

Proposition 2.13 is proved.

PROPOSITION 2.14. The asymptotic time complexity of method №1 is
\begin{spacing}{0.5}\end{spacing}\[O\left(m\cdot (\log _{2}^{6} m)\cdot \log _{2}^{O(1)} (\log _{2} m) \right).\]
\begin{spacing}{1.2}\end{spacing}

PROOF. Given Conjecture 2.1, let's say that
\begin{spacing}{0.5}\end{spacing}
 \[D(m^{4} )<p^{2} (j)<(m+1)^{4}.\]
\begin{spacing}{1.2}\end{spacing}

Let $t_{r} (m)$ -- number of actions in Step $r$ of method №1, $r=\overline{1;4}$;

{\rm \; \; \; \;}$t(m)$ -- time complexity of method №1 at input $m\in {\rm \mathbb{N}}: m\ge 3$.

\textbf{Step 1. } The number of digits (length) of $m$ equals $\lfloor\lg m  \rfloor+1$. 
\begin{spacing}{0.5}\end{spacing}
\begin{equation} \label{(21)} 
\lfloor\lg m\rfloor+1{\rm \;}\mathop{\le}\limits^{{\rm \eqref{(7)}}} {\rm \;}\lg m+1{\rm \;}\mathop{<}\limits^{m \ge 3}{\rm \;}\lg m+\lg(m^{3})=\lg m+3 \cdot \lg m=4 \cdot \lg m.
\end{equation} 
\begin{spacing}{1.2}\end{spacing}

Given (21), the multiplying of $m$ by itself requires not more than $O(\lg^{2}m)$ steps. All other arithmetic operations over intermediate results equire not more than $O(\lg m)$ steps. Therefore, the computing a number $\overline{n}$ of numbers of the form $6h\pm 1$ less than or equal to $m^{2}$ is going to take not more than $O(\lg^{2}m)$ steps. Thus, there may be
\begin{spacing}{0.5}\end{spacing}
\[t_{1} (m)=O(\lg^{2}m).\]
\begin{spacing}{1.2}\end{spacing}

\textbf{Step 2.} The asymptotic time complexity of Lenstra-Pomerance primality test is $\widetilde{O}(\log _{2}^{6} x)$ at input $x$ (see [4]), where 
\begin{spacing}{0.5}\end{spacing}
\[\widetilde{O}(y)=O\left(y\cdot (\log _{2} y)^{O(1)} \right).\]
\begin{spacing}{1.2}\end{spacing}

 The largest possible number $\Delta (m)$ of numbers $6h\pm 1$, which lie between $(m-1)^{2}$ and $m^{2}$, pass the primality test, equals $\nu(m^{2})-\nu((m-1)^{2})$. Then similar to \eqref{(19)}, we find an upper estimate for $\Delta (m)$:
\begin{spacing}{-0.2}\end{spacing}
\[\Delta (m)=\nu(m^{2})-\nu((m-1)^{2} )<\frac{m^{2} -(m-1)^{2} }{3} +2=\frac{(m+m-1)\cdot (m-m+1)}{3} +2=\] 
\begin{equation} \label{(22)} 
=\frac{2m-1}{3}+2=\frac{2m}{3} +\frac{5}{3} <\frac{2m}{3} +\frac{2m}{3} =\frac{4}{3} \cdot m{\rm \; \; }(m\ge 3){\rm \; \; }\Rightarrow {\rm \; \; }\Delta (m)<\frac{4}{3} \cdot m{\rm \; \; }(m\ge 3).  
\end{equation} 
\begin{spacing}{1.7}\end{spacing}

 Let 
\[f(\bar{n})=\max_{\scriptstyle \substack{(m-1)^{2} \le f(n)\le m^{2} \\ n\in {\rm \mathbb{N}}}}{\rm \;}f(n).\]

Introduce $t_{2,{\rm \; }i+1} (m)$ --  number of actions in Step 2 of method №1, when the number $f(\bar{n}-i), 0\le i\le \Delta (m)-1$, passes the primality test. Then, for $t_{2,{\rm \; }i+1} (m),{\rm \; }i=\overline{0;\Delta (m)-1},$ there are estimates:
\begin{spacing}{0.5}\end{spacing}
\[t_{2,{\rm \; }i+1} (m)=\widetilde{O}(\log _{2}^{6} f(\bar{n}-i))=O\left((\log _{2}^{6} f(\bar{n}-i))\cdot \log _{2}^{O(1)} (\log _{2}^{6} f(\bar{n}-i))\right).\]
\begin{spacing}{1.5}\end{spacing}

So, there exist constants $C_{1} >0,{\rm \; }C_{2} >0$ and $m_{0} \in {\rm \mathbb{N}} :{\rm \; }m_{0} \ge 3$, such that for all $m\ge m_{0} $:
\begin{spacing}{0.5}\end{spacing}
\[t_{2,{\rm \; }i+1} (m)\le C_{1} \cdot (\log _{2}^{6} f(\bar{n}-i))\cdot \log _{2}^{C_{2} } (\log _{2}^{6} f(\bar{n}-i)), i=\overline{0;\Delta (m)-1}.\] 
\[{\rm \; \;}C_{1} \cdot (\log _{2}^{6} f(\bar{n}-i))\cdot \log _{2}^{C_{2} } (\log _{2}^{6} f(\bar{n}-i)) \le C_{1} \cdot (\log _{2}^{6} (m^{2}))\cdot \log _{2}^{C_{2} } (\log _{2}^{6} (m^{2}))=\]
\[=C_{1} \cdot (2\cdot \log _{2} m)^{6} \cdot (6 \cdot \log _{2} (2\cdot \log _{2} m))^{{\rm \; }C_{2} }{\rm \; \;}\mathop{<}\limits^{m \ge 3}{\rm \; \;}C_{1} \cdot 2^{6} \cdot (\log _{2}^{6} m) \cdot (6 \cdot \log _{2} (\log_{2}^{3}m))^{{\rm \; }C_{2}}= \] 
\[=C_{1} \cdot 2^{6} \cdot (\log _{2}^{6} m) \cdot (3 \cdot \log _{2}(\log_{2}^{6}m))^{{\rm \; }C_{2}}=C_{1} \cdot 2^{6} \cdot 3^{C_{2}} \cdot (\log _{2}^{6} m) \cdot \log _{2}^{C_{2}} (\log_{2}^{6}m).\]

Introduce a constant $C_{3} \equiv C_{1} \cdot 2^{6} \cdot 3^{C_{2}} $. Then,
\begin{spacing}{0.5}\end{spacing}
\[(\forall m\ge m_{0} ) \left(t_{2,{\rm \; }i+1} (m)<C_{3}  \cdot (\log _{2}^{6} m) \cdot \log _{2}^{C_{2}} (\log_{2}^{6}m) \right), {\rm \;}i=\overline{0;\Delta (m)-1}.\] 

It follows that
\begin{spacing}{0.5}\end{spacing}
\[t_{2,{\rm \; }i+1} (m)=O\left((\log _{2}^{6} m) \cdot \log _{2}^{O(1)} (\log_{2}^{6}m)\right)=\widetilde{O}(\log _{2}^{6} m), {\rm \;}i=\overline{0;\Delta (m)-1}.\]
\begin{spacing}{1.2}\end{spacing}

Estimate the number of actions $t_{2} (m)$ in Step 2 of method №1 for all $m\ge m_{0}$.
\begin{spacing}{0.8}\end{spacing}
\[t_{2} (m)=\sum _{i=0}^{\Delta (m)-1}t_{2,i+1} (m)=\sum _{i=0}^{\Delta (m)-1}\widetilde{O}(\log _{2}^{6} m)=\widetilde{O}(\log _{2}^{6} m) \cdot \sum _{i=0}^{\Delta (m)-1}1=\] 
\begin{spacing}{0.4}\end{spacing}
\[=\widetilde{O}(\log _{2}^{6} m) \cdot \Delta (m) {\rm \;}\mathop{=}\limits^{{\rm \eqref{(22)}}}{\rm \;}\widetilde{O}(\log _{2}^{6} m) \cdot O(m)=O\left((\log _{2}^{6} m) \cdot \log _{2}^{O(1)} (\log_{2}^{6}m)\right) \cdot O(m)=\]
\[=O\left(m \cdot (\log _{2}^{6} m) \cdot \log _{2}^{O(1)} (\log_{2}^{6}m)\right)=O\left(O\left(m \cdot (\log _{2}^{6} m) \cdot \log _{2}^{O(1)} (\log_{2}m)\right)\right)=\]
\begin{spacing}{-0.2}\end{spacing}
\[=O\left(m \cdot (\log _{2}^{6} m) \cdot \log _{2}^{O(1)} (\log_{2}m)\right){\rm \; \; }\Rightarrow{\rm \; \; }t_{2} (m)=O\left(m \cdot (\log _{2}^{6} m) \cdot \log _{2}^{O(1)} (\log_{2}m)\right).\]

\textbf{Step 3.} The multiplying of $m$ by itself requires not more than $O(\lg^{2}m)$ steps. Then the multiplying of $m^2$ by itself also requires not more than $O(\lg^{2}m)$ steps, since 
\begin{spacing}{0.5}\end{spacing}
\[O(\lg^{2}(m^2))=O(4 \cdot \lg^{2}m)=O(\lg^{2}m).\]
\begin{spacing}{1.2}\end{spacing}

Further,
\[(\lfloor\lg(m^{4})\rfloor+1) \cdot (\lfloor\lg p(k_{1})\rfloor+1){\rm \;}\mathop{\le}\limits^{{\rm \eqref{(7)}}}{\rm \;}(\lg(m^{4})+1) \cdot (\lg p(k_{1})+1){\rm \;}\mathop{<}\limits^{{\rm \eqref{(14)}}}\]
\[< (\lg(m^{4})+1) \cdot (\lg(m^{2})+1){\rm \;}\mathop{<}\limits^{m \ge 3}{\rm \;} (\lg(m^{4})+\lg(m^{3})) \cdot (\lg(m^{2})+\lg(m^{3}))=\]
\[=(\lg (m^{7})) \cdot \lg (m^{5})=7 \cdot (\lg m) \cdot 5 \cdot \lg m=35 \cdot \lg^{2} m.\]

It follows that the division $m^4$ by $p(k_{1})$ with remainder is going to take not more than $O(\lg^{2} m)$ steps. On the basis of \eqref{(21)}, the asymptotic time complexity of computing the value of $C(m)$ is $O(\lg m)$.
\[\lg \left(\left\lfloor \frac{m^{4}}{p(k_{1})} \right\rfloor \right) {\rm \;}\mathop{\le}\limits^{{\rm \eqref{(7)}}}{\rm \;}\lg \left( \frac{m^{4}}{p(k_{1})} \right) {\rm \;}\mathop{<}\limits^{{\rm \eqref{(14)}}}{\rm \;} \lg \left( \frac{m^{4}}{(m-1)^{2}} \right) {\rm \;}\mathop{<}\limits^{m \ge 3}{\rm \;}\]
\begin{spacing}{0.2}\end{spacing} 
\[<\lg \left( \frac{m^{4}}{\left(m-\left(1-\frac{1}{\sqrt{3}} \right) \cdot m\right)^{2}} \right)=\lg \left( \frac{m^{4}}{\left(\frac{m}{\sqrt{3}}\right)^{2}} \right)=\lg(3 \cdot m^{2}) \le \lg(m^{3})=3 \cdot \lg m.\]
\begin{spacing}{1.5}\end{spacing} 
Therefore, the computing a number $n_{0}$ of $T$-matrix upper element $W(m^{4})$ of number $m^{4}$ (in $k_{1}$-row) with the computing $m^{4}$ and $\left\lfloor \frac{m^{4}}{p(k_{1})} \right\rfloor$ is going to take not more than $O(\lg^{2}m)$ steps. Thus, there may be
\begin{spacing}{0.5}\end{spacing}
\[t_{3} (m)=O(\lg^{2}m).\]
\begin{spacing}{1.2}\end{spacing}

\textbf{Step 4. }Let $t_{4,{\rm \; }i+1} (m)$ -- number of actions in Step 4 of method №1, when $f(n_{0} +i)$, where $0\le i\le \Delta n_{0}$, passes the primality test. Then, given the asymptotic time complexity of Lenstra-Pomerance primality test, for $t_{4,{\rm \; }i+1} (m),{\rm \; }i=\overline{0;\Delta n_{0} }$, there are estimates:
\begin{spacing}{-0.2}\end{spacing}
\[t_{4,{\rm \; }i+1} (m)=\widetilde{O}(\log _{2}^{6} f(n_{0} +i))=O\left((\log _{2}^{6} f(n_{0} +i))\cdot \log _{2}^{O(1)} (\log _{2}^{6}f(n_{0} +i))\right).\] 

So, there exist constants $C_{4}>0,{\rm \; }C_{5}>0$ and $m_{0} \in {\rm \mathbb{N}}: m_{0} \ge 3$, such that for all $m\ge m_{0}$:
\begin{spacing}{0.3}\end{spacing}
\[t_{4,{\rm \; }i+1} (m)\le C_{4} \cdot (\log _{2}^{6} f(n_{0} +i))\cdot \log _{2}^{C_{5} } (\log _{2}^{6} f(n_{0} +i)),{\rm \; }i=\overline{0;\Delta n_{0} }.\] 

There are the following inequalities:
\begin{spacing}{0.5}\end{spacing}
\begin{equation} \label{(23)}
m^{2}<f(n_{0} +i)\le p(j)<(m+1)^{2} ,{\rm \; }i=\overline{0;\Delta n_{0} }.
\end{equation}

Really,${\rm \; }m^{2} =\frac{m^{4} }{m^{2} } {\rm \; \; }\mathop{<}\limits^{\eqref{(14)}} {\rm \; }\frac{m^{4} }{p(k_{1} )} $, where by Theorem 2.12, $p(k_{1})\not {\rm |}{\rm \; }m^{4}$. So, $\frac{m^{4} }{p(k_{1} )} \notin {\rm \mathbb{N}}$.
\[\frac{m^{4} }{p(k_{1} )}<\frac{a(k_{1} ,n_{0} +i)}{p(k_{1} )}{\rm \;}\mathop{=}\limits^{\eqref{(2)}}{\rm \;}f(n_{0} +i){\rm \; ,\; }i=\overline{0;\Delta n_{0} }. \] 

Therefore, $m^{2} <f(n_{0} +i){\rm ,\; }i=\overline{0;\Delta n_{0} }$. Clear that $f(n_{0} +i)\le p(j),{\rm \; }i=\overline{0;\Delta n_{0} }$, and 
\begin{spacing}{0.5}\end{spacing}
\[f(n_{0} +i)=p(j){\rm \; \; }\Leftrightarrow {\rm \; \; }i=\Delta n_{0} .\] 

Since $(m+1)^{2} \notin {\rm \mathbb{P}}$, there is a strict inequality $p(j)<(m+1)^{2}$.
\[{\rm \; \;}C_{4} \cdot (\log _{2}^{6} f(n_{0} +i))\cdot \log _{2}^{C_{5} } (\log _{2}^{6} f(n_{0} +i)){\rm \; }\mathop{\le }\limits^{{\rm \eqref{(23)}}}{\rm \;} C_{4} \cdot (\log _{2}^{6} p(j))\cdot \log _{2}^{C_{5} } (\log _{2}^{6}p(j)){\rm \;}\mathop{<}\limits^{\eqref{(23)}} {\rm \;}\]
\begin{spacing}{-0.2}\end{spacing}
\[<C_{4} \cdot (\log _{2}^{6}((m+1)^{2}))\cdot \log _{2}^{C_{5} } (\log _{2}^{6}((m+1)^{2}))=\]
\[=C_{4} \cdot (2\cdot \log _{2} (m+1))^{6} \cdot (6 \cdot \log_{2} (2\cdot \log _{2} (m+1)))^{C_5}{\rm \;}\mathop{<}\limits^{m \ge 3}{\rm \;}\]
\[<C_{4} \cdot (2\cdot \log _{2} (m^{2}))^{6} \cdot (6 \cdot \log_{2} (2\cdot \log _{2} (m^{2})))^{C_5}=\]
\[=C_{4} \cdot 4^{6} \cdot (\log _{2}^{6} m) \cdot (6 \cdot \log_{2} (4\cdot \log _{2}m))^{C_5}{\rm \;}\mathop{<}\limits^{m \ge 3}{\rm \;}C_{4} \cdot 4^{6} \cdot (\log _{2}^{6} m) \cdot (6 \cdot \log_{2} (\log _{2}^{5}m))^{C_5}=\]
\[=C_{4} \cdot 4^{6} \cdot (\log _{2}^{6} m) \cdot (5 \cdot \log_{2} (\log _{2}^{6}m))^{C_5}=C_{4} \cdot 4^{6} \cdot 5^{C_5} \cdot (\log _{2}^{6} m) \cdot \log_{2}^{C_5} (\log _{2}^{6}m).\]

Introduce a constant $C_{6} \equiv C_{4} \cdot 4^{6} \cdot 5^{C_{5}}$. Then,
\begin{spacing}{0.5}\end{spacing}
\[(\forall m\ge m_{0} ) \left(t_{4,{\rm \; }i+1} (m)< C_{6}\cdot (\log _{2}^{6}m) \cdot \log _{2}^{C_{5}} (\log _{2}^{6}m)\right),{\rm \; }i=\overline{0;\Delta n_{0}}.\] 

It follows that 
\begin{spacing}{0.2}\end{spacing}
\[t_{4,{\rm \; }i+1} (m)=O\left((\log _{2}^{6}m) \cdot \log _{2}^{O(1)} (\log _{2}^{6} m)\right)=\widetilde{O}(\log _{2}^{6} m),{\rm \; }i=\overline{0;\Delta n_{0}}.\]
\[\sum _{i=0}^{\Delta n_{0} }t_{4,i+1} (m)=\sum _{i=0}^{\Delta n_{0}}\widetilde{O}(\log _{2}^{6} m)=\widetilde{O}(\log _{2}^{6} m) \cdot \sum _{i=0}^{\Delta n_{0}}1=\widetilde{O}(\log _{2}^{6} m) \cdot (\Delta n_{0} +1).\]

Note that on the basis of Lemma 2.1, 
\begin{spacing}{-0.5}\end{spacing}
\[f(n_{0}+i)\notin {\rm \mathbb{P}} \backslash \{ 2;3;5\} {\rm \;}\Leftrightarrow {\rm \;}a(k_{1};n_{0}+i)\notin D_{T}{\rm \;}\Leftrightarrow {\rm \;}a(k_{1};n_{0}+i)\in nD_{T},{\rm \;}i=\overline{0;\Delta n_{0}-1}.\]
\begin{spacing}{0.7}\end{spacing}

 Therefore, 
\begin{spacing}{0.5}\end{spacing}
\[\Delta n_{0} +1=\nu_{k_{1} } (D(m^{4} ))-\nu_{k_{1} } (m^{4} ).\] 

\begin{spacing}{0.5}\end{spacing}
\[\sum _{i=0}^{\Delta n_{0} }t_{4,i+1} (m)=\widetilde{O}(\log _{2}^{6} m) \cdot (\nu_{k_{1} } (D(m^{4} ))-\nu_{k_{1} }(m^{4} )){\rm \; }\mathop{=}\limits^{\eqref{(18)}}{\rm \; }\widetilde{O}(\log _{2}^{6} m) \cdot  O(m)=\]
\[=O\left((\log _{2}^{6}m) \cdot \log _{2}^{O(1)} (\log _{2}^{6}m)\right) \cdot O(m)=O\left(m \cdot (\log _{2}^{6}m) \cdot \log _{2}^{O(1)} (\log _{2}^{6}m)\right)=\]
\[=O\left(O\left(m \cdot (\log _{2}^{6} m) \cdot \log _{2}^{O(1)} (\log_{2}m)\right)\right)=O\left(m \cdot (\log _{2}^{6} m) \cdot \log _{2}^{O(1)} (\log_{2}m)\right){\rm \;}\Rightarrow\]
\begin{spacing}{0.5}\end{spacing}
\[\Rightarrow{\rm \; \;} \sum _{i=0}^{\Delta n_{0} }t_{4,i+1} (m)=O\left(m \cdot (\log _{2}^{6}m) \cdot \log _{2}^{O(1)} (\log _{2}m)\right).\]
\[(\lfloor \lg p(k_{1})\rfloor+1) \cdot (\lfloor\lg p(j)\rfloor+1){\rm \;}\mathop{\le}\limits^{{\rm \eqref{(7)}}}{\rm \;}(\lg p(k_{1})+1) \cdot (\lg p(j)+1)<\]
\begin{spacing}{0.5}\end{spacing}
\[<(\lg p(j)+1)^{2}{\rm \; \;}\mathop{<}\limits^{m \ge 3}{\rm \;}(2 \cdot \lg p(j))^{2}=4 \cdot \lg^{2} p(j){\rm \;}\mathop{<}\limits^{\eqref{(23)}} {\rm \;}4 \cdot \lg^{2}((m+1)^2)=\]
\begin{spacing}{0.5}\end{spacing}
\[=16 \cdot \lg^{2}(m+1){\rm \;}\mathop{<}\limits^{m \ge 3}{\rm \;}16 \cdot \lg^{2}(m^{2})=64 \cdot \lg^{2} m.\]

Therefore, the computing an element $D(m^{4})$ by multiplying of $p(k_{1}) {\rm \; }(p(k_{1})<p(j))$ by $p(j)$ is going to take not more than $O(\lg^{2}m)$ steps. Estimate the possible number of actions $t_{4} (m)$ in Step 4 of method №1 for all $m\ge m_{0}$.
\[t_{4} (m)=\sum _{i=0}^{\Delta n_{0} }t_{4,i+1} (m)+O(\lg^{2}m)=\]
\[{\rm \; \; \; \;}=O\left(m \cdot (\log _{2}^{6}m) \cdot \log _{2}^{O(1)} (\log _{2}m)\right)+O(\lg^{2}m)=O\left(m \cdot (\log _{2}^{6}m) \cdot \log _{2}^{O(1)} (\log _{2}m)\right){\rm \;}\Rightarrow\]
\[\Rightarrow{\rm \; \;}t_{4} (m)=O\left(m \cdot (\log _{2}^{6}m) \cdot \log _{2}^{O(1)} (\log _{2}m)\right).\]
\[t(m)=\sum _{i=1}^{4}t_{i} (m)=O(\lg^{2}m)+O \left(m\cdot (\log _{2}^{6} m)\cdot \log _{2}^{O(1)}(\log _{2}m) \right)+\]
\[+O(\lg^{2}m)+O \left(m\cdot (\log _{2}^{6} m)\cdot \log _{2}^{O(1)}(\log _{2}m) \right)=O \left(m\cdot (\log _{2}^{6} m)\cdot \log _{2}^{O(1)}(\log _{2}m) \right){\rm \;}\Rightarrow\]
\begin{spacing}{0.5}\end{spacing}
\[\Rightarrow{\rm \;}t(m)=O \left(m\cdot (\log _{2}^{6} m)\cdot \log _{2}^{O(1)}(\log _{2}m) \right).\] 

Proposition 2.14 is proved.

PROPOSITION 2.15. The asymptotic time complexity of finding numbers $k_{1}$, $j$ after all the steps of method №1 is $O(m^{1+o(1)})$ at input $m$.

PROOF. Initially, all steps of method №1 are completed. 

Let $t'_{1} (m)$ -- number of actions for finding a number $\overline{k}_{1}$;

$t'_{2} (m)$ -- number of actions for finding a number $\overline{j}$;

$t'_{3} (m)$ -- number of actions over numbers $\overline{k}_{1}$, $\overline{j}$;

$t'(m)$ -- time complexity of finding numbers $k_{1}$, $j$ at input $m$.
\begin{spacing}{1.7}\end{spacing}
The asymptotic time complexity of the Lagarias-Odlyzko analytical method is $O\left(x^{\frac{1}{2}+o(1)}\right)$ 
\begin{spacing}{1.7}\end{spacing}
\noindent (see [5], [6]). Then,
\begin{spacing}{0.5}\end{spacing}
\[t'_{1} (m)=O\left(p(k_{1})^{\frac{1}{2}+o(1)}\right) {\rm \; \; } \wedge {\rm \; \; } t'_{2} (m)=O\left(p(j)^{\frac{1}{2}+o(1)}\right).\] 
\begin{spacing}{1.5}\end{spacing}
$1){\rm \;} t'_{1} (m)=O\left(p(k_{1})^{\frac{1}{2}+o(1)}\right)$ means that there exist constants $C'_{1}>0$ and $m_{0} \in {\rm \mathbb{N}}: m_{0} \ge 3$, 
\begin{spacing}{1.8}\end{spacing}
\noindent such that for all $m\ge m_{0}$:
\begin{spacing}{0.5}\end{spacing}
\[t'_{1} (m) \le C'_{1} \cdot p(k_{1})^{\frac{1}{2}+o(1)}.\]
\begin{spacing}{-0.4}\end{spacing}
\[ C'_{1} \cdot p(k_{1})^{\frac{1}{2}+o(1)}{\rm \;} \mathop{<}\limits^{\eqref{(14)}}{\rm \;}C'_{1} \cdot (m^{2})^{\frac{1}{2}+o(1)}=C'_{1} \cdot m^{1+2 \cdot o(1)}=C'_{1} \cdot m^{1+o(1)},{\rm \;}m\ge m_{0}{\rm \; \; }\Rightarrow\]
\[\Rightarrow{\rm \; \; }t'_{1} (m)<C'_{1} \cdot m^{1+o(1)},{\rm \;}m\ge m_{0}{\rm \; \; }\Rightarrow{\rm \; \; }t'_{1} (m)=O\left(m^{1+o(1)}\right).\]
\begin{spacing}{1.5}\end{spacing}
$2){\rm \;} t'_{2} (m)=O\left(p(j)^{\frac{1}{2}+o(1)}\right)$ means that there exist constants $C'_{2}>0$ and $m_{0} \in {\rm \mathbb{N}}: m_{0} \ge 3$,
\begin{spacing}{1.8}\end{spacing} 
\noindent such that for all $m\ge m_{0}$: 
\begin{spacing}{0.5}\end{spacing} 
\[t'_{2} (m) \le C'_{2} \cdot p(j)^{\frac{1}{2}+o(1)}.\]
\begin{spacing}{-0.5}\end{spacing}
\[ C'_{2} \cdot p(j)^{\frac{1}{2}+o(1)}{\rm \;} \mathop{<}\limits^{\eqref{(23)}}{\rm \;}C'_{2} \cdot ((m+1)^{2})^{\frac{1}{2}+o(1)}=C'_{2} \cdot (m+1)^{1+2 \cdot o(1)}=C'_{2} \cdot (m+1)^{1+o(1)}<\]
\begin{spacing}{-0.2}\end{spacing}
\[<C'_{2} \cdot (1.34 \cdot m)^{1+o(1)}=C'_{2} \cdot 1.34^{1+o(1)} \cdot m^{1+o(1)}, {\rm \;}m\ge m_{0}{\rm \; \; }\Rightarrow\]
\[\Rightarrow{\rm \; \; }t'_{2} (m)<C'_{2} \cdot 1.34^{1+o(1)} \cdot m^{1+o(1)}, {\rm \;}m\ge m_{0}{\rm \;}\Rightarrow{\rm \;}t'_{2} (m)=O\left(1.34^{1+o(1)} \cdot m^{1+o(1)}\right){\rm \;}\Rightarrow{\rm \;}\]
\[\Rightarrow{\rm \;}t'_{2} (m)=O\left(1.34^{1+O(1)} \cdot m^{1+o(1)}\right).\]

The latter means that there exist constants $C'_{2}>0, C'_{3}>0${\rm \;} and {\rm \;}$m_{0} \in {\rm \mathbb{N}}: m_{0} \ge 3$,
such that for all $m\ge m_{0}$:
\begin{spacing}{0.2}\end{spacing}
\[t'_{2} (m)<C'_{2} \cdot 1.34^{1+C'_{3}} \cdot m^{1+o(1)}.\]

Introduce a constant $C'_{4} \equiv C'_{2} \cdot 1.34^{1+C'_{3}}$. Then,
\begin{spacing}{0.5}\end{spacing}
\[(\forall m\ge m_{0} )(t'_{2} (m)<C'_{4} \cdot m^{1+o(1)}).\]

It follows that
\[t'_{2} (m)=O\left(m^{1+o(1)}\right).\]
\begin{spacing}{1.2}\end{spacing}

$3) {\rm \;}\overline{k}_{1}<\overline {j}{\rm \; \;} \wedge {\rm \;} \lg \overline {j}=\lg \pi(p(j))<\lg p(j)<4 \cdot \lg m {\rm \; \;} \Rightarrow {\rm \; \;} t'_{3} (m)=O(\lg m)$.
\[{\rm \; \;}t'(m)=\sum _{i=1}^{3}t'_{i} (m)=O\left(m^{1+o(1)}\right)+O\left(m^{1+o(1)}\right)+O(\lg m)=O\left(m^{1+o(1)}\right) {\rm \; \; }\Rightarrow\]
\begin{spacing}{0.3}\end{spacing}
\[\Rightarrow{\rm \; \; }t'(m)=O\left(m^{1+o(1)}\right).\]

Proposition 2.15 is proved.

COROLLARY 2.16. The asymptotic time complexity of method №1 with the finding 
\begin{spacing}{1.7}\end{spacing}
\noindent numbers $k_{1}$, $j$ is $O\left(m\cdot (\log _{2}^{6} m)\cdot \log _{2}^{O(1)} (\log _{2} m) \right)$.
\begin{spacing}{2}\end{spacing}
PROOF. With the same notations in proof of Proposition 2.14 and proof of Proposition 2.15, 
\begin{spacing}{-0.3}\end{spacing}
\[t'(m)=O\left(m^{1+o(1)}\right)=O\left(m \cdot m^{o(1)}\right)=O\left(m \cdot O\left( (\log _{2}^{6} m)\cdot \log _{2}^{O(1)} (\log _{2} m) \right)\right)=\]
\[=O\left(O\left(m\cdot (\log _{2}^{6} m)\cdot \log _{2}^{O(1)} (\log _{2} m) \right)\right)=O\left(m\cdot (\log _{2}^{6} m)\cdot \log _{2}^{O(1)} (\log _{2} m) \right) {\rm \;}\Rightarrow\]
\begin{spacing}{0.8}\end{spacing}
\[\Rightarrow {\rm \;}t(m)+t'(m)=O\left(m\cdot (\log _{2}^{6} m)\cdot \log _{2}^{O(1)} (\log _{2} m) \right).\] 
\begin{spacing}{1.6}\end{spacing}

Corollary 2.16 is proved.

CONCLUSION 2.2. Method №1 with the finding numbers $k_{1}$, $j$ is a polynomial-time method.

EXAMPLE 2.1. Find a $T$-matrix upper defining element $D(10^{4})$ of number $10^{4}$. 

Find numbers $k_{1}$, $j$.

SOLUTION.\textbf{ Input:} $m=10$.

\textbf{Step 1}. Compute a number $\overline{n}$ of numbers of the form $6h\pm 1$ less than or equal to $m^{2}=100$.
\[\overline{n}\equiv \nu(m^{2})=\left\lfloor \frac{10^{2}+2}{3} \right\rfloor -\left\lfloor \frac{10^{2}{\rm \; }\% 6}{4} \right\rfloor +\left\lfloor \frac{10^{2}{\rm \; }\% 6}{5} \right\rfloor -1=34-1+0-1=32.\] 

\textbf{Step 2. }Using the Lenstra-Pomerance primality test,
\begin{spacing}{1}\end{spacing}
\[f(\overline{n})=3\cdot 32+\frac{3-(-1)^{32} }{2} =97\in {\rm \mathbb{P}} \text{ when } i=0.\]

In that case, $\Delta \overline{n}=0$. Therefore, $p(k_{1})=f(\overline{n}-\Delta \overline{n})=97$. 

Within $T$-matrix,
\begin{spacing}{0.3}\end{spacing}
\[{\rm \;}9^{4}<p^{2} (k_{1})=9409<10^{4}.\]

\textbf{Step 3. }Compute a number $n_{0} $ of $T$-matrix upper element $W(m^{4})$ of $m^{4}$ (in $k_{1}$-row).

\[n_{0} =C\left(\left\lfloor \frac{m^{4} }{p(k_{1} )} \right\rfloor \right)=C(103)=\left\lfloor \frac{105}{3} \right\rfloor -\left\lfloor \frac{103{\rm \; }\%{\rm \; } 6}{4} \right\rfloor +\left\lfloor \frac{103{\rm \; }\%{\rm \; } 6}{5} \right\rfloor =35.\] 

\textbf{Step 4. }Using the Lenstra-Pomerance primality test, we find an element $D(m^{4})$ of $T$-matrix.
\[i=0:{\rm \; \; }f(n_{0} )=3\cdot 35+\frac{3-(-1)^{35} }{2} =105+2=107\in {\rm \mathbb{P}} {\rm \;}\Rightarrow {\rm \;}\Delta n_{0} =0.\] 
\[D(m^{4} )=p(k_{1} )\cdot f(n_{0} +\Delta n_{0} )=97\cdot 107=10379<(m+1)^{4} =11^{4} =14641.\] 

The obtained prime number $p(j)=107$ lies between $10^{2} $ and $11^{2}$. Within $T$-matrix,
\begin{spacing}{0.5}\end{spacing}
\[10^{4}<p^{2} (j)=11449<11^{4}.\]
\begin{spacing}{1.2}\end{spacing}

\textbf{Output:} $D(m^{4})=10379$. 

Using the Lagarias-Odlyzko analytical method for prime numbers {\rm \;}$p(k_{1} )=97, p(j)=107$,
\begin{spacing}{0.5}\end{spacing}
\[\overline{k_{1} }=\pi (p(k_{1} ))=25{\rm \;}\Rightarrow {\rm \;}k_{1} =23.\] 
\[\overline{j}=\pi (p(j))=28{\rm \;}\Rightarrow {\rm \;}j =26.\] 

\begin{center}
\textbf{Table 1. }Fragment of $T$-matrix for Example 2.1
\end{center}
{\rm \; \; \;}\includegraphics*[scale=0.64]{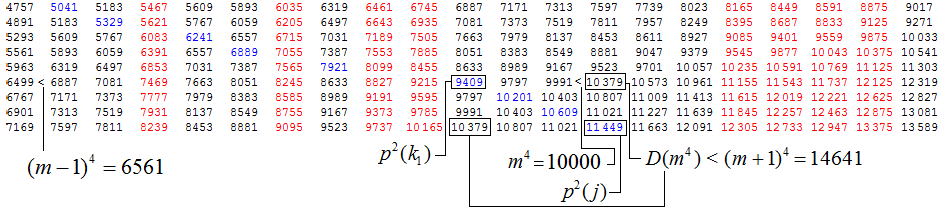}

Example 2.1 is considered.

Now we'll deal with the question of choice of $T$-matrix leading element $p^{2} (k_{1} )$ on the basis of Step 2 of method №1. Let any number $m\in {\rm \mathbb{N}}, {\rm \; }m\ge 3$ is chosen. From some number $k>1$ select all $d\in {\rm \mathbb{N}}$ leading elements $p^{2} (k+i{\rm )}$, $i=\overline{0;d-1}$ that lie between $(m-1)^{4}$ and $m^{4}$. Using Definition 2.3, we find the upper defining elements $D_{k+i} (m^{4})$ of number $m^{4}$, and
\begin{spacing}{0.5}\end{spacing}
\[D_{k+i} (m^{4})<(m+1)^{4} ,{\rm \; }i=\overline{0;d-1}.\] 
\begin{spacing}{1.1}\end{spacing}
Since  $(m-1)^{4} ,{\rm \; }m^{4} \notin {\rm M} _{T} $, there are strict inequalities:
\begin{spacing}{0.3}\end{spacing}
\[(m-1)^{4} <p^{2} (k+i)<m^{4} <D_{k+i} (m^{4} ),{\rm \; \; }i=\overline{0;d-1}.\] 

Use Theorem 2.2 for elements $D_{k+i} (m^{4})$. As a result, there is a «transition down» of element $D_{k+i} (m^{4})$ from  ($k+i$)-row $(i=\overline{0;d-1})$ to the appropriate $j_{i}$-row (some of them may coincide) of $T$-matrix:
\begin{spacing}{0.5}\end{spacing}
\[k+i<j_{i} {\rm \; \;}\wedge {\rm \; \;}D_{k+i} (m^{4} )<p^{2} (j_{i} ){\rm \; \;}\wedge {\rm \; \;}a(j_{i} ;\# _{k+i} (p^{2} (k+i)))=D_{k+i} (m^{4} ){\rm \; \;}\wedge \] 
\[\wedge {\rm \; \;}a(j_{i} ;\# _{k+i} (D_{k+i} (m^{4} )))=p^{2} (j_{i} ),{\rm \; \; }i=\overline{0;d-1}.\]

 In turn, given Conjecture 2.1, there may be the following 3 situations.

\textbf{Situation №1. }$d=1$. Then there exists only 1 leading element $p^{2}(k)$ between $(m-1)^{4}$ and $m^{4}$. When the «transition down» happens from $k$-row to $j_{0} $-row, there are inequalities:
\begin{spacing}{0.5}\end{spacing}
\[m^{4} <D_{k} (m^{4} )<p^{2} (j_{0} )<(m+1)^{4} .\]  

 \textbf{Situation №2.} $d\ne 1{\rm \; \; \; }\wedge {\rm \; \; \; }m^{4} <D_{k+i} (m^{4} )<(m+1)^{4} <p^{2} (j_{i} ),{\rm \; }i=\overline{0;h}{\rm ,\; }h<d-1{\rm \; \; \; }\wedge $
\begin{spacing}{0.5}\end{spacing}
\[\wedge {\rm \; \; \; }m^{4} <D_{k+i} (m^{4} )<p^{2} (j_{i} )<(m+1)^{4} ,{\rm \; }i=\overline{h+1;d-1}.\] 

 \textbf{Situation №3. }$d\ne 1{\rm \; \; \; }\wedge {\rm \; \; \; }m^{4} <D_{k+i} (m^{4} )<(m+1)^{4} <p^{2} (j_{i} ),{\rm \; }i=\overline{0;d-2}{\rm \; \; \; }\wedge $
\begin{spacing}{0.5}\end{spacing}
\[\wedge {\rm \; \; \; }m^{4} <D_{k+d-1} (m^{4} )<p^{2} (j_{d-1} )<(m+1)^{4} .\] 

Really, if there is a chain of inequalities 
\begin{spacing}{0.5}\end{spacing}
\[m^{4} <D_{k+i} (m^{4} )<(m+1)^{4} <p^{2} (j_{i} ),{\rm \; \; }i=\overline{0;d-1},\]
then the result may be that the prime numbers don't exist between $m^{2} $ and $(m+1)^{2}$. This would mean that the Legendre's conjecture is false. That's why Definition 2.1 of $T$-matrix upper defining element of number was introduced. Hence question of choice of $T$-matrix leading element $p^{2}(k_{1})$ is solved. 

The benefit of choosing this leading element is that $k_{1}$-row of $T$-matrix contains all defining elements $a(k_{1};n_{i})>p^{2}(k_{1}),{\rm \; }i=\overline{1;q_{m}},{\rm \; } n_{1}<n_{2}<...<n_{q_{m}}$; for each of them with appropriate «transition» to other $j_{i}$-row of $T$-matrix one of the chains of inequalities is done:
\begin{equation} \label{(24)} 
\left[\begin{array}{l} {a(k_{1};n_{i})<m^{4} <p^{2} (j_{i})<(m+1)^{4}} \\ {m^{4} <a(k_{1} ;n_{i})<p^{2} (j_{i})<(m+1)^{4}} \end{array}\right. \text{.} 
\end{equation} 

Consequently, we will consider the next paragraph.
\begin{spacing}{0.2}\end{spacing}
\begin{center}
\textbf{\large 3. «Active» set and «critical» element for numbers $(m-1)^{4}, m^{4}$ $ (m\ge 3)$}
\end{center}
\begin{spacing}{0.78}\end{spacing}

Let any natural number $m\ge 3$ is chosen, and the leading element $p^{2} (k_{1})$ is found on the basis of \eqref{(14)}. Assume that Conjecture 2.1 is true.

 Let us consider an ordered set $(D_{T_{k_{1} } } ;\le )$, where
\begin{spacing}{0.5}\end{spacing}
\[D_{T_{k_{1} } } \equiv \{ p(k)\cdot p(k_{1} ),{\rm \; }k\ge 2{\rm \} }.\]
\begin{spacing}{1.2}\end{spacing}

Since all elements of $(D_{T_{k_{1} } } ;\le )$ are pairwise comparable, $(D_{T_{k_{1} } } ;\le )$ is a linearly ordered set, and an appropriate relation $\le $ is a relation of linear order.

The following principles are known.

THEOREM 3.1 (greatest element principle). Every nonempty finite set of natural numbers has a greatest element.

THEOREM 3.2 (least element principle). Every nonempty set of natural numbers has a least element.

PROPERTY 3.1. 1){\rm \;}$7\cdot p(k_{1})$ is the least element of set $D_{T_{k_{1}}}$.

                        2) The set $D_{T_{k_{1}}} $ has not the greatest element.

                        3) $(D_{T_{k_{1}}};\le)$ is well-ordered set.

PROOF. 1) $D_{T_{k_{1} } } \subset {\rm \mathbb{N}} {\rm \; \; }\wedge {\rm \; \; }D_{T_{k_{1} } } \ne \varnothing $. By Theorem 3.2, the set $D_{T_{k_{1}}} $ has the least element. 

In turn,
\[(\forall x\in D_{T_{k_{1}}})(p(2)\cdot p(k_{1})\le x), \text{ where } p(2)=7.\]
As such, Property 3.1, 1) is true.

2) Introduce a linearly ordered set $({\rm \mathbb{P}} \backslash \{2;3;5\} ;\le )$. By Euclid's theorem (see [3]), the set $\mathbb{P}$ is infinite. So, the set $\mathbb{P}$ has not the greatest element. Then the set ${\rm \mathbb{P}} \backslash \{2;3;5\} $ has not the greatest element. It follows that the set $D_{T_{k_{1}}}$ has not the greatest element.

3) $D_{T_{k_{1} } } \subset {\rm \mathbb{N}} {\rm \; \; }\wedge {\rm \; \; }D_{T_{k_{1} } } \ne \varnothing $.  Then every subset of set $D_{T_{k_{1} }}$ is also a subset of set ${\rm \mathbb{N}}$. Then, by Theorem 3.2, every nonempty subset of set $D_{T_{k_{1}}} $ has the least element. That means that $(D_{T_{k_{1} } } ;\le )$ is well-ordered set. 

Property 3.1 is proved.

DEFINITION 3.1. A set  ${\rm H} _{(m-1)^{4} ,{\rm \; }m^{4} }\subset$  $D_{T_{k_{1} } } $ of all defining elements  $a(k_{1};n_{i})>p^{2} (k_{1})$ (from $k_{1}$-row of $T$-matrix), $i=\overline{1;q_{m}}$ , which satisfy \eqref{(24)}, is called an «active» set for numbers $(m-1)^{4}$, $m^{4}$ ($m\ge 3$).

Let 
\begin{spacing}{-0.1}\end{spacing}
\begin{equation} \label{(25)} 
{\rm H} _{(m-1)^{4} ,{\rm \; }m^{4} } \equiv {\rm \; \{ }a(k_{1} ;n_{1} ),{\rm .}..,{\rm \; }a(k_{1} ;n_{q_{m}} )\} ,\text{ where } a(k_{1} ;n_{1} )<...<a(k_{1} ;n_{q_{m}} ).
\end{equation} 
\begin{equation} \label{(26)} 
\text{GCD}({\rm H} _{(m-1)^{4} ,{\rm \; }m^{4}})\equiv{ \text{GCD}\rm (}a(k_{1} ;n_{1} ),...,a(k_{1} ;n_{q_{m}})).                                         
\end{equation} 

DEFINITION 3.2. A defining element  $C_{(m-1)^{4} ,{\rm \; }m^{4}} \equiv a(k_{1} ;n_{q_{m}+1})\notin {\rm H}_{(m-1)^{4} ,{\rm \; }m^{4} } $, next to a defining element $a(k_{1} ;n_{q_{m}} )\in {\rm H} _{(m-1)^{4} ,{\rm \; }m^{4} }$, is called a «critical» element for numbers $(m-1)^{4} ,{\rm \; }m^{4}$ ($m\ge 3$).

DEFINITION 3.3. A «transition» of the defining element $a(k_{1};n_{i})>p^{2} (k_{1})$ from $k_{1}$-row to $j_{i}$-row ($j_{i}>k_{1}$) of $T$-matrix with some $i\in {\rm \mathbb{N}}$ is called successful if $a(k_{1};n_{i})\in {\rm H} _{(m-1)^{4} ,{\rm \; }m^{4} } $. Otherwise, that is if $a(k_{1} ;n_{i})\notin {\rm H}_{(m-1)^{4} ,{\rm \; }m^{4}}$, this «transition» is called unsuccessful.

PROPERTY 3.2. The «active» set ${\rm H} _{(m-1)^{4} ,{\rm \; }m^{4} }$ is finite.

PROOF. The number $q_{m}$ of prime numbers between $m^{2} $ and $(m+1)^{2} $ is finite. Then number of defining elements which satisfy \eqref{(24)} and greater than the leading element $p^{2} (k_{1})$ is also finite. So, the set ${\rm H} _{(m-1)^{4} ,{\rm \; }m^{4}}$ is finite. 

Property 3.2 is proved.

Note that by Theorem 3.1 and Theorem 3.2, the set ${\rm H} _{(m-1)^{4} ,{\rm \; }m^{4} }$ has the greatest and least elements. It follows from \eqref{(25)} that 
\begin{spacing}{0.2}\end{spacing}
\[\min {\rm H} _{(m-1)^{4} ,{\rm \; }m^{4} } =a(k_{1} ;n_{1} ){\rm \; \;}\wedge {\rm \; \;}\max {\rm H} _{(m-1)^{4} ,{\rm \; }m^{4} } =a(k_{1} ;n_{q_{m}}).\]

PROPERTY 3.3.
\[\min {\rm H} _{(m-1)^{4} ,{\rm \; }m^{4} }=D(p^{2} (k_{1})).\]

This equality follows from Definition 3.1 and Definition 2.1.

PROPOSITION 3.3.
\begin{equation} \label{(27)} 
\text{GCD}({\rm H}_{(m-1)^{4},{\rm \; }m^{4}})=p(k_{1}).
\end{equation} 

PROOF. Definition 3.1 makes it clear that all elements of set ${\rm H} _{(m-1)^{4} ,{\rm \; }m^{4} }$ are defining and 
\begin{spacing}{0.5}\end{spacing}
\[a(k_{1};n_{i})>p^{2} (k_{1}), i=\overline{1;q_{m}}.\] 
\begin{spacing}{1.2}\end{spacing}

Then we can use Theorem 2.2 for them. As a result,
\begin{spacing}{0.5}\end{spacing}
\[a(k_{1} ;n_{i}){\rm \; \; }\mathop{=}\limits^{\eqref{(5)}}{\rm \; \; }a(j_{i};\# _{k_{1}} (p^{2} (k_{1}))){\rm \; \; }\mathop{=}\limits^{\eqref{(2)}}{\rm \; \; }p(j_{i}) \cdot f(\# _{k_{1}} (p^{2} (k_{1})))=\]
\begin{equation} \label{(28)} 
=p(j_{i}) \cdot p(k_{1})= p(k_{1}) \cdot p(j_{i})\Rightarrow {\rm \; \;}a(k_{1} ;n_{i})=p(k_{1} )\cdot p(j_{i}),{\rm \; }i=\overline{1;q_{m}}.  
\end{equation} 
\begin{spacing}{0.2}\end{spacing}
\[a(k_{1} ;n_{1} )<...<a(k_{1} ;n_{q_{m}}){\rm \; \; }\mathop{\Leftrightarrow }\limits^{\eqref{(28)}} {\rm \; \; }p(k_{1} )\cdot p(j_{1} )<...<p(k_{1} )\cdot p(j_{q_{m}}){\rm \; \; }\mathop{\Leftrightarrow }\limits^{p(k_{1} )>0}\] 
\[\Leftrightarrow{\rm \; \;} p(j_{1} )<...<p(j_{q_{m}} ){\rm \; \; }\Rightarrow {\rm \; \; }\text{GCD}(p(j_{1} ),...,p(j_{q_{m}} ))=1{\rm \; \; }\Leftrightarrow\]
\[\Leftrightarrow {\rm \; \; }p(k_{1} )\cdot \text{GCD}(p(j_{1} ),...,p(j_{q_{m}} ))=p(k_{1} ).\] 
\[p(k_{1} )\cdot \text{GCD}(p(j_{1} ),...,p(j_{q_{m}}))= \text{GCD}(p(k_{1} )\cdot p(j_{1} ),...,p(k_{1} )\cdot p(j_{q_{m}})){\rm \; \; }\mathop{=}\limits^{(28)} \] 
\[=\text{GCD}(a(k_{1} ;n_{1} ),...,a(k_{1} ;n_{q_{m}}))\mathop{=}\limits^{\eqref{(26)}}  \text{GCD}({\rm H}_{(m-1)^{4} ,{\rm \; }m^{4}}){\rm \; \;}\Rightarrow {\rm \; \;}\text{GCD}({\rm H}_{(m-1)^{4},{\rm \; }m^{4}})=p(k_{1}).\] 
\begin{spacing}{1}\end{spacing}

Proposition 3.3 is proved.

PROPOSITION 3.4. If divide all elements of «active» set  ${\rm H} _{(m-1)^{4} ,{\rm \; }m^{4} }$ by $\text{GCD}({\rm H} _{(m-1)^{4} ,{\rm \; }m^{4}})$, then get all the different prime numbers that lie between $m^{2}$ and $(m+1)^{2}$ $(m\ge 3$).

PROOF. By Definition 3.1, every element of «active» set ${\rm H} _{(m-1)^{4} ,{\rm \; }m^{4} }$ satisfies \eqref{(24)}. Therefore, for each of them with appropriate «transition» from $k_{1}$-row to other $j_{i}$-row of $T$-matrix:
\begin{spacing}{0.5}\end{spacing}
\[m^{4} <p^{2} (j_{i} )<(m+1)^{4} ,{\rm \; }i=\overline{1;q_{m}}{\rm \; \; \; }\Leftrightarrow {\rm \; \; \; }m^{2} <p(j_{i} )<(m+1)^{2} ,{\rm \; }i=\overline{1;q_{m}}. \]

Using the beginning of the proof of Proposition 3.3, we come to presentations \eqref{(28)} of the defining elements $a(k_{1} ;n_{i}),{\rm \; }i=\overline{1;q_{m}}$. Express the prime numbers $p(j_{i})$:
\begin{spacing}{0.8}\end{spacing}
 \begin{equation} \label{(29)} 
p(j_{i})=\frac{a(k_{1} ;n_{i} )}{p(k_{1})} ,{\rm \; }i=\overline{1;q_{m}}{\rm \; \; }\mathop{\Leftrightarrow }\limits^{\eqref{(27)}} {\rm \; \; }p(j_{i})=\frac{a(k_{1} ;n_{i} )}{\text{GCD}({\rm H}_{(m-1)^{4},{\rm \; }m^{4}})},{\rm \; }i=\overline{1;q_{m}}.
\end{equation} 
\begin{spacing}{0.5}\end{spacing}
\begin{equation} \label{(30)} 
a(k_{1} ;n_{1} )<...<a(k_{1} ;n_{q_{m}} ){\rm \; \; }\mathop{\Leftrightarrow }\limits^{\scriptstyle \substack { \eqref{(28)}, \\  p(k_{1})>0}}{\rm \; \; } p(j_{1})<...<p(j_{q_{m}}).  
\end{equation} 

Thus, given \eqref{(29)} and \eqref{(30)}, we make sure that Proposition 3.4 is true. 

Proposition 3.4 is proved.

PROPOSITION 3.5. The $T$-matrix upper defining element $D(m^{4})<(m+1)^{4}$ of number $m^{4}$ $(m\ge 3)$ belongs to the «active» set ${\rm H} _{(m-1)^{4} ,{\rm \; }m^{4} }$ for numbers $(m-1)^{4} ,{\rm \; }m^{4}$.

This proposition is true, if Conjecture 2.1 is true. We'll show it.
\begin{spacing}{1}\end{spacing}
\[\left(\forall m\in {\rm \mathbb{N}}:{\rm \;} m\ge 3\right)\left(m^{2} <\frac{D(m^{4} )}{p(k_{1} )} <(m+1)^{2} \right){\rm \; \;}\Leftrightarrow\]

\[\Leftrightarrow {\rm \; \;} \left(\forall m\in {\rm \mathbb{N}}:{\rm \;} m\ge 3\right)\left(m^{4} <\frac{D^{2}(m^{4} )}{p^{2}(k_{1} )} <(m+1)^{4}\right).\]

It follows from Definition 2.1 and Definition 1.1 that  
\begin{spacing}{1}\end{spacing}
\[\frac{D^{2}(m^{4})}{p^{2}(k_{1})}=p^{2}(j).\]
\[\frac{D^{2}(m^{4})}{p^{2}(k_{1})}{\rm \; \; }\mathop{>}\limits^{\eqref{(14)}} {\rm \; \; }\frac{D^{2}(m^{4})}{m^{4}}>\frac{D^{2}(m^{4})}{D(m^{4})}=D(m^{4}){\rm \; \;}\Rightarrow {\rm \; \;}D(m^{4})<p^{2}(j).\]

Using Theorem 2.6 for defining element $D(m^{4})<p^{2}(j)$, we make sure that all elements of «active» set ${\rm H} _{(m-1)^{4} ,{\rm \; }m^{4} }$ will be located in $k_{1}$-row $(k_{1}<j)$ of $T$-matrix. One such element is $D(m^{4})$. Really, when the element $D(m^{4})$ «moves down» from $k_{1} $-row to $j$-row ($j>k_{1}$) of $T$-matrix, there is a chain of inequalities:
\begin{spacing}{0.3}\end{spacing}
\begin{equation} \label{(31)} 
m^{4} <D(m^{4} )<p^{2} (j)<(m+1)^{4} .  
\end{equation} 
\begin{spacing}{1}\end{spacing}

Therefore, condition \eqref{(24)} holds in relation to the element $D(m^{4} )>p^{2} (k_{1} )$. By Definition 3.1, that means that $D(m^{4} )\in {\rm H} _{(m-1)^{4} ,{\rm \; }m^{4} }$. 

Thus, by Definition 3.3, the «transition down» of element $D(m^{4})$ from $k_{1} $-row to $j$-row ($j>k_{1}$) of $T$-matrix will be successful.

COROLLARY 3.6. The $T$-matrix upper defining element $D(m^{4})<(m+1)^{4}$ of number $m^{4} $ $(m\ge 3)$ is not «critical» for numbers $(m-1)^{4} ,{\rm \; }m^{4} $.

PROPOSITION 3.7. The «transition» of «critical» element $C_{(m-1)^{4} ,{\rm \; }m^{4}}$ for numbers 

\noindent $(m-1)^{4} ,{\rm \; }m^{4}$ $(m\ge 3)$ from $k_{1}$-row to $j_{q_{m}+1}$-row ($j_{q_{m}+1}>k_{1}$) of $T$-matrix is unsuccessful.

This proposition follows from Definition 3.2 and Definition 3.3.

PROPOSITION 3.8. The defining elements $a(k_{1} ;n_{i}){\rm ,\; }i=\overline{1;s_{m}},{\rm \; }s_{m}<q_{m}$, lying between the leading element $p^{2} (k_{1})$ and the $T$-matrix upper defining element $D(m^{4})\equiv a(k_{1} ;n_{s_{m}+1} )$ $(a(k_{1} ;n_{s_{m}+1} )<(m+1)^{4})$ of number $m^{4}$ $(m\ge 3)$, are elements of «active» set ${\rm H} _{(m-1)^{4},{\rm \; }m^{4} } $.

PROOF. Let ${\rm H} '_{(m-1)^{4} ,{\rm \; }m^{4} }$ is a set of all defining elements $a(k_{1} ;n_{i})$ such that
\begin{spacing}{0.7}\end{spacing}
\[p^{2} (k_{1})<a(k_{1};n_{i})<D(m^{4}){\rm ,\; }i=\overline{1;s_{m}},{\rm \; }s_{m}<q_{m}.\]
\begin{spacing}{1.2}\end{spacing}

Using Theorem 2.2 for each element of set ${\rm H} '_{(m-1)^{4} ,{\rm \; }m^{4} }$, we get
\begin{spacing}{0.5}\end{spacing}
\[k_{1} <j_{i} {\rm \; \;}\wedge {\rm \; \;}a(k_{1} ;n_{i} )<p^{2} (j_{i} ){\rm \; \; }\wedge {\rm \; \; }a(j_{i} ;{\rm \; }\# _{k_{1} } (p^{2} (k_{1} )))=a(k_{1} ;n_{i} ){\rm \; \; }\wedge\]
\[\wedge{\rm \; \;}a(j_{i} ;n_{i} )=p^{2} (j_{i} ){\rm ,\; }i=\overline{1;s_{m}}.\] 
\begin{spacing}{1.2}\end{spacing}

As opposed to the element $D(m^{4})$ for each of them: 
\begin{spacing}{0.3}\end{spacing}
\[a(k_{1} ;n_{i} )<m^{4} ,{\rm \; }i=\overline{1;s_{m}}.\] 
\begin{spacing}{1.2}\end{spacing}
Now we need to show that inequalities $m^{4} <p^{2} (j_{i} )<(m+1)^{4} ,{\rm \; }i=\overline{1;s_{m}}$; are true. 

Assume the converse. Then consider 2 cases.

\textbf{Case 1.} $p^{2} (j_{r})<m^{4}$ with some $r\in {\rm \mathbb{N}}:{\rm \; }1\le r\le s_{m}$. 

Within $T$-matrix,
\begin{spacing}{0.5}\end{spacing}
\begin{equation} \label{(32)} 
p(k_{1})=\max_{\scriptstyle \substack{(m-1)^{2}<p(k)<m^{2} \\ k>1}}{\rm \; \; } p(k){\rm \; \; }\Leftrightarrow{\rm \; \; }p^{2}(k_{1})=\max_{\scriptstyle \substack{(m-1)^{4}<p^{2}(k)<m^{4} \\  k>1}}{\rm \; \; }p^{2}(k).
\end{equation} 
\begin{spacing}{1.5}\end{spacing}

It follows from Property 1.1 and inequality $k_{1} <j_{r}$ that $p^{2} (k_{1})<p^{2} (j_{r} )$. Then,

\[(m-1)^{4} <p^{2} (k_{1})<p^{2} (j_{r})<m^{4} .\] 
\begin{spacing}{1.2}\end{spacing}

As a result, a contradiction to the maximality of the leading element $p^{2} (k_{1})$ from \eqref{(32)}.

\textbf{Case 2. }$(m+1)^{4}<p^{2} (j_{r})$ with some $r\in {\rm \mathbb{N}} :{\rm \; }1\le r\le s_{m}$. 

The $T$-matrix upper defining element $D(m^{4})$ of number $m^{4}$ satisfies inequalities \eqref{(31)}. 

Then, 
\begin{spacing}{0.2}\end{spacing}
\[m^{4}<D(m^{4})<p^{2} (j)<(m+1)^{4}<p^{2} (j_{r} ){\rm \; \; }\Rightarrow{\rm \; \; } p^{2} (j)<p^{2} (j_{r} ). \] 
\begin{spacing}{1.2}\end{spacing}

From condition of Proposition 3.8,
\begin{spacing}{0.3}\end{spacing}
\[p^{2} (k_{1})<a(k_{1} ;n_{r})<D(m^{4}).\]
\begin{spacing}{1.2}\end{spacing}

Then it follows from Theorem 2.2 for defining elements $a(k_{1};n_{r})$ and $D(m^{4})$ that $j_{r}<j$. So, by Property 1.1, $p^{2} (j_{r})<p^{2}(j)$. As a result, a contradiction.

Thus,
\begin{spacing}{0.1}\end{spacing}
\[a(k_{1} ;n_{i} )<m^{4} <p^{2} (j_{i} )<(m+1)^{4} ,{\rm \; }i=\overline{1;s_{m}}.\] 
\begin{spacing}{1.2}\end{spacing}

Therefore, condition \eqref{(24)} holds in relation to every element $a(k_{1} ;n_{i})$ $(i=\overline{1;s_{m}})$. By Definition 3.1, that means that
\begin{spacing}{0.3}\end{spacing}
\[a(k_{1} ;n_{i} )\in {\rm H} _{(m-1)^{4} ,{\rm \; }m^{4} } ,{\rm \; }i=\overline{1;s_{m}}.\] 

Eventually,
\begin{spacing}{0.1}\end{spacing}
\[{\rm H} '_{(m-1)^{4} ,{\rm \; }m^{4} } \subset {\rm H} _{(m-1)^{4} ,{\rm \; }m^{4} }.\]

 Proposition 3.8 is proved.

EXAMPLE 3.1. Construct an «active» set ${\rm H} _{5^{4},{\rm \; }6^{4}}$, find a «critical» element $C_{5^{4},{\rm \; }6^{4}}$ for numbers $5^{4}, 6^{4}$.

SOLUTION. It follows from condition of Example 3.1 that $m=6$. Then,
\begin{spacing}{0.5}\end{spacing}
\[(m-1)^{4} =5^{4} =625,{\rm \; \; }m^{4} =6^{4} =1296,{\rm \; \; }(m+1)^{4} =7^{4} =2401.\] 

Using method №1 with the finding a number $k_{1}$ or the presentation of $T$-matrix, we compute:

1) $T$-matrix leading element $p^{2}(k_{1})$ satisfying \eqref{(32)}:
\begin{spacing}{0.5}\end{spacing}
\[p(k_{1} )=31{\rm \; \; }\Leftrightarrow{\rm \; \; } p^{2} (k_{1} )=961.\]

2) number of $k_{1}$-row of $T$-matrix:
\begin{spacing}{0.1}\end{spacing}
\[k_{1}=9.\]
\begin{spacing}{1}\end{spacing}

3) $T$-matrix upper defining element $D(m^{4} )<(m+1)^{4}$ of number $m^{4}$:
\begin{spacing}{0.4}\end{spacing}
\[D(m^{4} )=1333.\]

First, we construct an «active» set ${\rm H} _{(m-1)^{4} ,{\rm \; }m^{4} }$ for numbers $(m-1)^{4} ,{\rm \; }m^{4}$. 

1) By Proposition 3.5,
\[D(m^{4} )=1333\in {\rm H} _{(m-1)^{4} ,{\rm \; }m^{4} }.\]
\begin{spacing}{1.2}\end{spacing}

2) By Proposition 3.8, the defining elements $a(k_{1} ;n_{1} )=1147,{\rm \; }a(k_{1} ;n_{2} )=1271$, lying between the leading element $p^{2} (k_{1} )$ and the element $D(m^{4})$, are elements of set ${\rm H} _{(m-1)^{4} ,{\rm \; }m^{4} }$. In this case, let 
\begin{spacing}{0.3}\end{spacing}
\[a(k_{1} ;n_{3})\equiv D(m^{4}).\]

3) The defining element $a(k_{1};n_{4})=1457$ immediately follows after the defining element $a(k_{1} ;n_{3})=1333$ in $k_{1}$-row of $T$-matrix. It follows from Definition 1.1 and Theorem 2.2 that 
\begin{spacing}{0.8}\end{spacing}
\begin{equation} \label{(33)} 
p(j_{i} )=\frac{a(k_{1} ;n_{i} )}{p(k_{1} )} ,{\rm \; }i=\overline{1;q_{m}+1}.  
\end{equation}
\begin{spacing}{1.1}\end{spacing}
\[p(j_{4} ){\rm \; \; }\mathop{=}\limits^{\eqref{(33)}} {\rm \; \; }\frac{1457}{31} =47<(m+1)^{2} =7^{2} =49.\] 
\begin{spacing}{1.5}\end{spacing}
As a result, condition \eqref{(24)} holds for element $a(k_{1};n_{4})$. So, by Definition 3.1,
\begin{spacing}{0.5}\end{spacing}
\[a(k_{1} ;n_{4} )\in {\rm H} _{(m-1)^{4} ,{\rm \; }m^{4} }.\]
\begin{spacing}{1.2}\end{spacing}

4) The not defining element 1519 immediately follows after the defining element $a(k_{1} ;n_{4})$ in $k_{1}$-row of $T$-matrix.

5) The defining element $a(k_{1};n_{5})=1643$ immediately follows after the not defining element 1519 in $k_{1}$-row of $T$-matrix.
\begin{spacing}{0.8}\end{spacing}
\[p(j_{5} ){\rm \; \; }\mathop{=}\limits^{\eqref{(33)}} {\rm \; \; }\frac{a(k_{1} ;n_{5} )}{p(k_{1} )} =\frac{1643}{31} =53>(m+1)^{2} .\] 

As a result, condition \eqref{(24)} does not hold for element $a(k_{1};n_{5})$. So, by Definition 3.1, 
\begin{spacing}{0.4}\end{spacing}
\[a(k_{1} ;n_{5} )\notin {\rm H} _{(m-1)^{4} ,{\rm \; }m^{4} }.\]

$a(k_{1};n_{5})$ is the next defining element after the defining element $a(k_{1} ;n_{4} )\in {\rm H} _{(m-1)^{4} ,{\rm \; }m^{4}}$. Then, by Definition 3.2, 
\begin{spacing}{0.5}\end{spacing}
\[C_{(m-1)^{4},{\rm \; }m^{4}}=a(k_{1} ;n_{5}).\]
\begin{spacing}{1.2}\end{spacing}

 Eventually,
\begin{spacing}{0.5}\end{spacing}
\[{\rm H} _{(m-1)^{4} ,{\rm \; }m^{4} } =\{ 1147;{\rm \; }1271;{\rm \; 1333;\; 1457}\} {\rm \; \; (}s_{m}={\rm 2,\; }q_{m}=4{\rm )},\]
\[{\rm \; } C_{(m-1)^{4} ,{\rm \; }m^{4} } =1643.\]

 Example 3.1 is considered.

Next, we’ll show the illustration of Example 3.1.

\textbf{}

\begin{center}
\textbf{Table 2.} Illustration of Example 3.1. 

«Active» set ${\rm H} _{5^{4} ,{\rm \; }6^{4}} $ for numbers $5^{4} ,{\rm \; }6^{4}$; «critical» element $C_{5^{4} ,{\rm \; }6^{4} }$ for numbers $5^{4},{\rm \; }6^{4}$.
\begin{spacing}{2.0}\end{spacing}
\includegraphics*[scale=0.7]{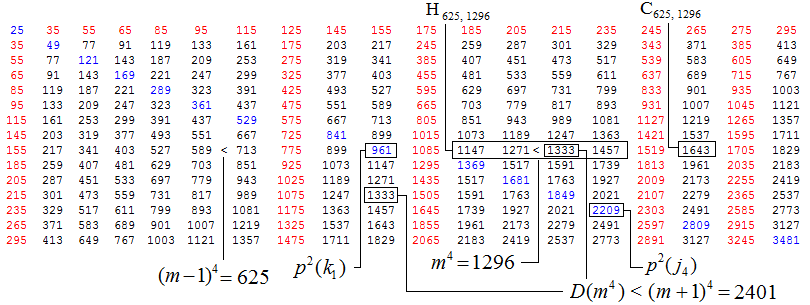}     
\end{center}
\begin{spacing}{1}\end{spacing}

 Now, let's look at an important example, when the $T$-matrix lower defining element $d(m^{4})$ of number $m^{4}$ coincides with the leading element $p^{2} (k_{1})$.

EXAMPLE 3.2. Construct an «active» set ${\rm H} _{4^{4},{\rm \; }5^{4}}$, find a «critical» element $C_{4^{4},{\rm \; }5^{4}}$ for numbers $4^{4}, 5^{4}$.

SOLUTION. It follows from condition of Example 3.2 that $m=5$. Then, 
\begin{spacing}{0.5}\end{spacing}
\[(m-1)^{4} =4^{4} =256,{\rm \; \; }m^{4} =5^{4} =625,{\rm \; \; }(m+1)^{4} =6^{4} =1296.\]

Using method №1 with the finding a number $k_{1}$ or the presentation of $T$-matrix, we compute:

1) $T$-matrix leading element $p^{2}(k_{1})$ satisfying \eqref{(32)}:
\begin{spacing}{0.5}\end{spacing}
\[p(k_{1} )=23{\rm \; \; }\Leftrightarrow{\rm \; \; } p^{2} (k_{1} )=529.\]
\begin{spacing}{1.2}\end{spacing}

2) number of $k_{1}$-row of $T$-matrix: 
\[k_{1} =7.\]
\begin{spacing}{1}\end{spacing}

3) $T$-matrix upper defining element $D(m^{4} )<(m+1)^{4}$ of number $m^{4}$:
\begin{spacing}{0.5}\end{spacing}
\[D(m^{4})=667.\]

First, we construct an «active» set ${\rm H} _{(m-1)^{4} ,{\rm \; }m^{4} }$ for numbers $(m-1)^{4} ,{\rm \; }m^{4}$. 

1) By Proposition 3.5, 
\[D(m^{4} )=667\in {\rm H} _{(m-1)^{4} ,{\rm \; }m^{4}}.\] 

2) As Table 3 shows, the defining elements don't exist between the leading element $p^{2} (k_{1})$ and the element $D(m^{4})$ in $k_{1}$-row of $T$-matrix, there exists only the not defining element 575. Therefore, by Definition 2.2, $d(m^{4})=p^{2} (k_{1})$. As such, let
\begin{spacing}{0.5}\end{spacing}
\[a(k_{1};n_{1})\equiv D(m^{4}).\]

3) The defining element $a(k_{1};n_{2})=713$ immediately follows after the defining element 

\noindent $a(k_{1};n_{1})=667$ in $k_{1}$-row of $T$-matrix.
\begin{spacing}{0.9}\end{spacing}
\[p(j_{2} ){\rm \; \; }\mathop{=}\limits^{\eqref{(33)}} {\rm \; \; }\frac{a(k_{1} ;n_{2} )}{p(k_{1} )} =\frac{713}{23} =31<(m+1)^{2} =6^{2} =36.\] 

As a result, condition \eqref{(24)} holds for element $a(k_{1} ;n_{2})$. So, by Definition 3.1,
\begin{spacing}{0.5}\end{spacing}
\[a(k_{1} ;n_{2})\in {\rm H} _{(m-1)^{4} ,{\rm \; }m^{4} }.\]
\begin{spacing}{1.2}\end{spacing}

4) The not defining element 805 immediately follows after the defining element $a(k_{1};n_{2})$ in $k_{1}$-row of $T$-matrix.

5) The defining element $a(k_{1} ;n_{3})=851$ immediately follows after the not defining element 805 in $k_{1}$-row of $T$-matrix.
\begin{spacing}{0.7}\end{spacing}
\[p(j_{3} ){\rm \; \; }\mathop{=}\limits^{\eqref{(33)}} {\rm \; \; }\frac{a(k_{1} ;n_{3} )}{p(k_{1} )} =\frac{851}{23} =37>(m+1)^{2} .\] 

As a result, condition \eqref{(24)} does not hold for element $a(k_{1};n_{3})$. So, by Definition 3.1, 
\begin{spacing}{0.5}\end{spacing}
\[a(k_{1} ;n_{3} )\notin {\rm H} _{(m-1)^{4} ,{\rm \; }m^{4} }.\]
\begin{spacing}{1.2}\end{spacing}

$a(k_{1};n_{3})$ is the next defining element after the defining element $a(k_{1} ;n_{2} )\in {\rm H} _{(m-1)^{4} ,{\rm \; }m^{4}}$. Then, by Definition 3.2, 
\begin{spacing}{0.2}\end{spacing}
\[C_{(m-1)^{4} ,{\rm \; }m^{4}}=a(k_{1};n_{3}).\]

Eventually,
\begin{spacing}{0.1}\end{spacing}
\[{\rm H} _{(m-1)^{4} ,{\rm \; }m^{4} } =\{ 667;{\rm \; }713\} {\rm \; (}s_{m}={\rm 0,\; }q_{m}=2{\rm )},{\rm \; }C_{(m-1)^{4} ,{\rm \; }m^{4} }=851.\]

Example 3.2 is considered.

Next, we’ll show the illustration of Example 3.2.
\begin{center}
\textbf{Table 3.} Illustration of Example 3.2.

«Active» set ${\rm H} _{4^{4} ,{\rm \; }5^{4}}$ for numbers $4^{4} ,{\rm \; }5^{4}$; «critical» element $C_{4^{4} ,{\rm \; }5^{4}}$ for numbers $4^{4} ,{\rm \; }5^{4}$.
\begin{spacing}{2.0}\end{spacing}
\includegraphics*[scale=0.7]{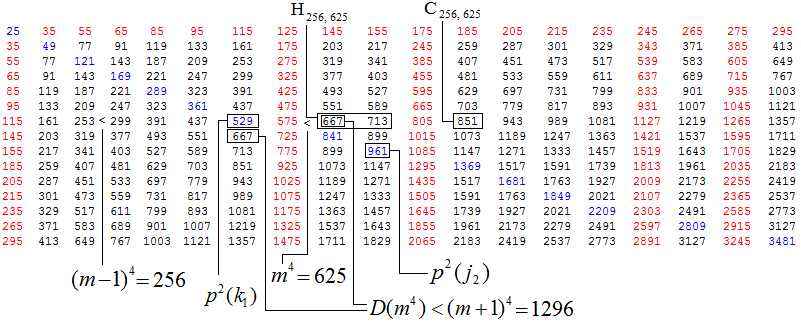}     
\end{center}
\begin{spacing}{0.7}\end{spacing}

Note that the equality $d(m^{4})=p^{2}(k_{1})$ does not affect the existence of prime number between $m^{2}$ and $(m+1)^{2}$, since in this case $D(p^{2} (k_{1}))=D(m^{4}) \in {\rm H} _{(m-1)^{4} ,{\rm \; }m^{4}}.$

THEOREM 3.9. 
\begin{spacing}{-0.5}\end{spacing}
\[\left(\forall x\in {\rm \mathbb{R}}\right)\left(x \ge \sqrt{3}{\rm \; \; }\Rightarrow {\rm \; \; }\pi ((x+1)^{2} )-\pi (x^{2} )=\pi _{{\rm M} _{T} } ((x+1)^{4} )-\pi _{{\rm M} _{T} } (x^{4} )\right) (\text{see } [1]).\]
\begin{spacing}{0.8}\end{spacing}

THEOREM 3.10 (about cardinality of «active» set ${\rm H} _{(m-1)^{4} ,{\rm \; }m^{4}}$). 
\begin{spacing}{0.5}\end{spacing}
\[\left(\forall m\ge 3\right)\left(\left|{\rm H} _{(m-1)^{4} ,{\rm \; }m^{4} } \right|=\pi _{{\rm M} _{T} } ((m+1)^{4} )-\pi _{{\rm M} _{T} } (m^{4} )\right).\]

PROOF. The cases $m=1$ and $m=2$ are not considered, since
\begin{spacing}{0.5}\end{spacing}
\[m\in \{ 1;2\} {\rm \; \; }\Rightarrow {\rm \; \; }(m-1)^{4} <m^{4}<p^{2} (1)={\rm 25}.\] 

It follows from Proposition 3.4 that 
\begin{spacing}{0.5}\end{spacing}
\[\left|{\rm H} _{(m-1)^{4} ,{\rm \; }m^{4}} \right|=q_{m}=\pi ((m+1)^{2})-\pi (m^{2}).\]

It follows from Theorem 3.9 that 
\begin{spacing}{0.5}\end{spacing}
\[\left(\forall m\in {\rm \mathbb{N}} \right)\left(m\ge 3{\rm \; \; }\Rightarrow {\rm \; \; }\pi ((m+1)^{2} )-\pi (m^{2} )=\pi _{{\rm M} _{T} } ((m+1)^{4} )-\pi _{{\rm M} _{T} } (m^{4} )\right).\]
\begin{spacing}{1.2}\end{spacing}

Theorem 3.10 is proved.

PROPOSITION 3.11. 1)${\rm \;}D(m^{4})<(m+1)^{4}{\rm \; \; }\Rightarrow {\rm \; \; }D(p^{2} (k_{1} ))\le D(m^{4} )$.

2) The element $D(p^{2} (k_{1}))$ is not «critical» for numbers $(m-1)^{4},{\rm \; }m^{4} $ $(m\ge 3)$.

PROOF. It follows from Definition 2.1 and \eqref{(32)} that 
\begin{spacing}{0.3}\end{spacing}
\[p^{2} (k_{1})<D(p^{2} (k_{1})){\rm \; \;}\wedge {\rm \; \;}p^{2} (k_{1})<m^{4}<D(m^{4}).\] 
\begin{spacing}{1.2}\end{spacing}

\textbf{Case 1. } There exists the defining element, lying between the leading element $p^{2}(k_{1})$ and the number $m^{4}$ in $k_{1}$-row of $T$-matrix. In this case, 
\begin{spacing}{0.3}\end{spacing}
\[p^{2} (k_{1} )<D(p^{2} (k_{1} ))<m^{4} <D(m^{4} ){\rm \; \; }\Rightarrow {\rm \; \; }D(p^{2} (k_{1} ))<D(m^{4} ).\]
\begin{spacing}{1}\end{spacing}
By Proposition 3.8,
\begin{spacing}{0.2}\end{spacing}
\[p^{2} (k_{1} )<D(p^{2} (k_{1} ))<D(m^{4} )<(m+1)^{4}{\rm \; \; }\Rightarrow {\rm \; \; }D(p^{2} (k_{1} ))\in {\rm H} _{(m-1)^{4} ,{\rm \; }m^{4} } .\] 

By Definition 3.2, that means that the element $D(p^{2} (k_{1}))$ is not «critical» for $(m-1)^{4} ,{\rm \; }m^{4} $.

\textbf{Case 2. } There doesn't exist the defining elements, lying between the leading element $p^{2}(k_{1})$ and $m^{4}$ in $k_{1}$-row of $T$-matrix. In this case, $D(p^{2} (k_{1}))=D(m^{4})$. From Corollary 3.6 we get that the element $D(p^{2} (k_{1} ))$ $(D(p^{2} (k_{1} ))<(m+1)^{4})$ is not «critical» for numbers $(m-1)^{4} ,{\rm \; }m^{4}$.

Proposition 3.11 is proved.

PROPOSITION 3.12. For $m \ge 3$, if Legendre's conjecture is true, then
\begin{spacing}{0.7}\end{spacing}
\[\min_{\scriptstyle \substack {m^2<p<(m+1)^2  \\  p \in \mathbb{P}}}{\rm \;} p=\frac{D(p^{2} (k_{1} ))}{p(k_{1})}.\]

PROOF. Use Property 3.3.
\begin{spacing}{0.9}\end{spacing}
\[D(p^{2} (k_{1} ))=\min{\rm H} _{(m-1)^{4} ,{\rm \; }m^{4} }{\rm \; \; }\Leftrightarrow {\rm \; \; }\frac{D(p^{2} (k_{1} ))}{p(k_{1})}=\frac{\min {\rm H} _{(m-1)^{4} ,{\rm \; }m^{4} }}{p(k_{1})}.\] 

\[\frac{\min {\rm H} _{(m-1)^{4} ,{\rm \; }m^{4} }}{p(k_{1})}{\rm \; \; }\mathop{=}\limits^{\eqref{(25)}} {\rm \; \; }\frac{a(k_{1} ;n_{1})}{p(k_{1})}{\rm \; \; }\mathop{=}\limits^{\eqref{(28)}} {\rm \; \; }p(j_{1}){\rm \; \; }\mathop{=}\limits^{\eqref{(30)}} \min_{\scriptstyle \substack { m^2<p<(m+1)^2  \\  p \in \mathbb{P}}}{\rm \;} p.\]

Proposition 3.12 is proved.

CONCLUSION 3.1.  For $m \ge 3$, if Legendre's conjecture is true, then
\begin{spacing}{0.8}\end{spacing}
\[\min_{\scriptstyle \substack{m^2<p<(m+1)^2  \\ p \in \mathbb{P}}}{\rm \;} p=\frac{\min{\rm H} _{(m-1)^{4} ,{\rm \; }m^{4} }}{\text{GCD}({\rm H}_{(m-1)^{4},{\rm \; }m^{4}})}.\]

\begin{center}
\textbf{\large 4. Major findings and conjectures}
\end{center}
\begin{spacing}{0.7}\end{spacing}

CONCLUSION 4.1. For $m\ge 2$, Legendre's conjecture is true ${\rm \; }\Leftrightarrow{\rm \;}$

\begin{center}
$(\exists q \in \mathbb{M}_{T})(q \in \left(m^{4}; (m+1)^{4})\right).$
\end{center}

CONCLUSION 4.2. For $m\ge 3$, Legendre's conjecture is true ${\rm \; }\Leftrightarrow{\rm \;}  {\rm H} _{(m-1)^{4} ,{\rm \; }m^{4} } \ne \varnothing$.

«WEAK» CONJECTURE 4.1. $D(p^{2} (k_{1}))\in {\rm H} _{(m-1)^{4} ,{\rm \; }m^{4}} $.

«STRONG» CONJECTURE 4.2. $D(m^{4} )\in {\rm H} _{(m-1)^{4} ,{\rm \; }m^{4} } $.

CONCLUSION 4.3. 1) Conjecture 2.1 is true ${\rm \; }\Leftrightarrow{\rm \;}$ «Strong» Conjecture 4.2 is true.

2) For $m\ge 3$, Legendre's conjecture is true ${\rm \; }\Leftrightarrow{\rm \;}$ «Weak» Conjecture 4.1 is true.

3)  For $m\ge 3$, «Weak» Conjecture 4.1 is true ${\rm \; }\Leftrightarrow{\rm \; \;} m^{2} <\frac{D(p^{2}(k_{1} ) )}{p(k_{1} )} <(m+1)^{2}.$

4) For $m\ge 3$, «Strong» Conjecture 4.2 is true ${\rm \; }\Rightarrow{\rm \;}$ Legendre's conjecture is true.

5) «Strong» Conjecture 4.2 is true ${\rm \; }\Rightarrow{\rm \;}$ «Weak» Conjecture 4.1 is true (this follows from Conclusion 4.3, 2) and Conclusion 4.3, 4)).

CONCLUSION 4.4. Only one of three outcomes of conjectures is true.

\textbf{Outcome №1. }Legendre's conjecture is true. ${\rm H} _{(m-1)^{4} ,{\rm \; }m^{4} } \ne \varnothing $.There exists the defining element, lying between the leading element $p^{2}(k_{1})$ and the number $m^{4}$ in $k_{1}$-row of $T$-matrix. The «transition» of element $D(m^{4})$ from $k_{1}$-row to $j$-row ($j>k_{1}$) of $T$-matrix is unsuccessful, the «transition» of element $D(p^{2}(k_{1}))$ from $k_{1}$-row to $j_{1}$-row ($j>j_{1}>k_{1}$) of $T$-matrix is successful («Strong» Conjecture 4.2 is false).

\textbf{Outcome №2. }Legendre's conjecture is false. ${\rm H} _{(m-1)^{4} ,{\rm \; }m^{4} }=\varnothing $. There doesn't exist the defining elements, lying between the leading element $p^{2}(k_{1})$ and $m^{4}$ in $k_{1}$-row of $T$-matrix. The «transition» of element $D(m^{4})$, where $D(m^{4})=D(p^{2}(k_{1}))$, from $k_{1}$-row to $j$-row ($j>k_{1}$) of $T$-matrix is unsuccessful («Strong» Conjecture 4.2 is false).

\textbf{Outcome №3. }Legendre's conjecture is true. ${\rm H} _{(m-1)^{4} ,{\rm \; }m^{4} } \ne \varnothing $.The «transition» of element $D(m^{4})$ from $k_{1}$-row to $j$-row ($j>k_{1}$) of $T$-matrix is successful, the «transition» of element $D(p^{2} (k_{1}))$ from $k_{1}$-row to $j_{1}$-row ($j\ge j_{1}>k_{1}$) of $T$-matrix is successful («Strong» Conjecture 4.2 is true).

\begin{center}
\textbf{\large 5. Conclusion}
\end{center}
\begin{spacing}{0.7}\end{spacing}

Theorem of presentation and basic definitions for $T$-matrix elements are formulated. New types of $T$-matrix elements, the most important of which is $T$-matrix upper defining element $D(b)$ of some real number $b \ge 49$, are introduced. Theorems and consequences from them related to the «transition» of defining elements, in particular, of element $D(b)$, from one row of $T$-matrix to another are proved. Way to go from the leading elements $p^{2}(k)$ and $p^{2}(k+1)$ to element $D_{k}(p^{2}(k)) {\rm \;} (k>1)$ is shown. Relation between the functions $\nu_{k}$ and $\nu$ is established. Formula for calculating the values of $\nu (m)$ for all $m \in \mathbb{N}$ is got. 

Method to compute the $T$-matrix upper defining element $D(m^{4})$ of number $m^{4}$ $(m\ge 3)$ with the finding numbers $k_{1}$, $j$ is developed. Asymptotic time complexity of this method is found. It has been shown that this method has polynomial running time. The problem of finding a prime number $p(j)$ between $m^{2}$ and $(m+1)^{2}$ is considered.

Conjectures about the ratios $\frac{D(m^4)}{p(k_{1})}$ and $\frac{D(p^{2}(k_{1}))}{p(k_{1})}$, lying between $m^{2}$ and $(m+1)^{2}{\rm \;}(m \ge 3)$, are proposed. The indivisibility of number $m^{4}{\rm \;}(m \ge 3)$ by prime number $p(k_{1})$ is proved.

Properties of sets $D_{T_{k_{1} } },{\rm \;}{\rm H} _{(m-1)^{4},{\rm \; }m^{4}}$ are explored. Propositions about «active» set and «critical» element for numbers $(m-1)^{4}$ and $m^{4}$ are proved, assuming that Conjecture 2.1 is true. The theorem about cardinality of «active» set ${\rm H} _{(m-1)^{4},{\rm \; }m^{4}}$ is proved. Important in finding elements of «active» set ${\rm H} _{(m-1)^{4},{\rm \; }m^{4}}$ is the inequality 

\begin{spacing}{0.5}\end{spacing}
\[D(p^{2} (k_{1} ))\le D(m^{4}) \text{, where }D(m^{4})<(m+1)^{4}.\]
\begin{spacing}{1.2}\end{spacing}

Two formulas for calculating the minimal prime number between $m^{2}$ and $(m+1)^{2}$ $(m\ge 3)$ are found, assuming that Legendre's conjecture is true.

«Weak» conjecture, «Strong» conjecture and their equivalent forms are got. Outcomes of these conjectures are described in connection with Legendre's conjecture.

\begin{center}
\textbf{\large References}
\end{center}
\begin{enumerate}
\item Garipov I. About one matrix of composite numbers and her applications //arXiv preprint arXiv:2012.15745. – 2020.

\item Семенов И. Л. Антье и мантисса. Сборник задач с решениями / Под ред. Е. В. Хорошиловой. М.: ИПМ им. М. В. Келдыша, 2015. -- 432 с.

\item Бухштаб А.А. Теория чисел, М., Просвещение, 1966.

\item H. W. Lenstra jr., C. Pomerance, ''Primality testing with Gaussian periods'', ver. of April 12, 2011.

\item Р. Крэндалл, К.Померанс. Простые числа: Криптографические и вычислительные аспекты. М.: УРСС, 2011.

\item J. C. Lagarias, A. M. Odlyzko, Computing $\pi (x)$: an analytic method, Journal of Algorithms 8, 1987, no. 2, 173-191. 

\item В. И. Зенкин. Распределение простых чисел. Элементарные методы, 2012. -- 112 c.
\end{enumerate}
\end{document}